\documentclass{article}[12pt]
\usepackage{amssymb}
\usepackage{amsmath}
\usepackage{amscd}

\title{Les vari\'et\'es sur le corps \`a un \'el\'ement}
\author{Christophe Soul\'e}

\def\ra{\rightarrow}
\def\longra{\longrightarrow}
\def\mpo{\mapsto}

\def\ify{\infty}
\def\ot{\otimes}
\def\sbs{\subset}
\def\sm{\simeq}
\def\ts{\times}

\def\a{\alpha}
\def\D{\Delta}
\def\G{\Gamma}
\def\lb{\lambda}
\def\Lb{\Lambda}
\def\Om{\Omega}
\def\s{\sigma}
\def\Si{\Sigma}
\def\t{\theta}
\def\ve{\varepsilon}
\def\vp{\varphi}
\def\z{\zeta}

\def\build#1_#2^#3{\mathrel{
\mathop{\kern 0pt#1}\limits_{#2}^{#3}}}

\begin{document}

\maketitle

Une fantaisie r\'ecurrente de plusieurs math\'ematiciens
(\cite{T}, \cite{Ma}, \cite{S}, \cite{KS},~$\ldots$)
 est l'existence d'un
``corps \`a un \'el\'ement'', not\'e ${\mathbb F}_1$, et d'une
g\'eom\'etrie alg\'ebrique sur ce corps. On pense par exemple
que le groupe des points de ${\rm SL}_N$ dans ${\mathbb F}_1$
est le groupe sym\'etrique des permutations de $N$ lettres, et
que ces $N$ lettres sont les points dans ${\mathbb F}_1$ de
l'espace projectif ${\mathbb P}^N$. Et l'on s'est aper\c cu
depuis longtemps que des formules connues pour les points d'un
groupe de Chevalley dans le corps fini ${\mathbb F}_q$, $q > 1$,
donnent par la sp\'ecialisation $q=1$ des formules vraies pour
le groupe de Weyl correspondant. Par ailleurs, l'analogie entre
corps de nombres et corps de fonctions incite \`a chercher un
corps de base pour la courbe affine ${\rm Spec} ({\mathbb Z})$.
Enfin, il y a un int\'er\^et croissant en g\'eom\'etrie
arithm\'etique pour les vari\'et\'es alg\'ebriques sur ${\mathbb
Z}$ issues de constructions combinatoires sur les ensembles
finis.

Le but de cet article est de proposer une d\'efinition des
vari\'et\'es sur ${\mathbb F}_1$. Pour ce faire, on part de
l'id\'ee qu'une vari\'et\'e $X$ (de type fini) sur ${\mathbb
F}_1$ doit avoir une extension des scalaires \`a ${\mathbb Z}$,
qui sera un sch\'ema $X_{{\mathbb Z}}$ de type fini sur
${\mathbb Z}$. Les points de $X$ (dans un anneau ad hoc) sont
alors une partie, qu'on supposera finie, de l'ensemble des
points de $X_{{\mathbb Z}}$. Par exemple, l'ensemble des points
du groupe multiplicatif ${\mathbb G}_m$ sur ${\mathbb F}_1$ dans
un anneau $R$ fini et plat sur ${\mathbb Z}$ est l'ensemble
(fini) des racines de l'unit\'e contenues dans $R$.

La vari\'et\'e $X_{{\mathbb Z}}$ doit \^etre enti\`erement
d\'etermin\'ee par $X$ (\`a isomorphisme unique pr\`es).
Autrement dit, les vari\'et\'es sur ${\mathbb Z}$ obtenues par
extension des scalaires de ${\mathbb F}_1$ \`a ${\mathbb Z}$
poss\`edent une description ``combinatoire
finie''. Notre r\'esultat principal (Th\'eor\`eme~1) est que les
vari\'et\'es toriques lisses peuvent \^etre d\'efinies sur
${\mathbb F}_1$. Un autre exemple, inspir\'e de la th\'eorie
d'Arakelov, consiste \`a associer \`a tout r\'eseau $\Lb \sm
{\mathbb Z}^d$ muni d'une norme hermitienne sur $\Lb
\build\ot_{{\mathbb Z}}^{} {\mathbb C}$ une vari\'et\'e affine
sur le corps ${\mathbb F}_1$. Cependant, nous n'avons pas su,
\`a ce stade, d\'efinir sur ${\mathbb F}_1$ les groupes de
Chevalley et les vari\'et\'es de drapeaux.

Pour en revenir \`a la d\'efinition d'une vari\'et\'e $X$ sur
${\mathbb F}_1$, il se trouve que la donn\'ee des points de $X$
ne suffit pas \`a d\'eterminer la vari\'et\'e
alg\'ebrique $X_{{\mathbb Z}}$ sur ${\mathbb Z}$. C'est pourquoi
nous sommes conduits \`a adjoindre \`a la donn\'ee des points
celle de fonctions, \`a savoir une ${\mathbb C}$-alg\`ebre
{\LARGE $a$}$_X$, qui sera dans nos exemples une alg\`ebre
de Banach commutative. La d\'efinition de $X$ 
permet aussi d'\'evaluer les
fonctions de {\LARGE $a$}$_X$ aux points de $X$. Dans le cas de
${\mathbb G}_m$ par exemple, l'alg\`ebre {\LARGE $a$}$_{{\mathbb
G}_m}$ est celle des fonctions continues sur le cercle unit\'e.

Une vari\'et\'e $V$ sur ${\mathbb Z}$ d\'efinit un ``objet sur
${\mathbb F}_1$'', constitu\'e du foncteur des points de $V$ et
de l'alg\`ebre des fonctions alg\'ebriques sur $V ({\mathbb
C})$. Si $X$ est une vari\'et\'e sur ${\mathbb F}_1$, son
extension $X_{{\mathbb Z}}$ \`a ${\mathbb Z}$ est un objet
initial parmi toutes les vari\'et\'es $V$ sur ${\mathbb Z}$
telles qu'il existe un morphisme d'objets sur ${\mathbb F}_1$ de
$X$ vers $V$ (D\'efinitions 3 et 5). L'existence de $X_{{\mathbb
Z}}$ est une propri\'et\'e non triviale de l'objet $X$.

Ce texte est organis\'e comme suit. Le premier paragraphe est un
bref historique du corps ${\mathbb F}_1$ et des sp\'eculations
auxquelles il a donn\'e lieu. Le second explique quels
raisonnements ont conduit aux d\'efinitions pr\'esent\'ees ici,
qui font l'objet du troisi\`eme paragraphe. La section~4 donne
quelques propri\'et\'es des vari\'et\'es sur ${\mathbb F}_1$,
qui indiquent une certaine coh\'erence de la th\'eorie
envisag\'ee. On peut par exemple d\'efinir des vari\'et\'es sur
${\mathbb F}_1$ par recollement (Proposition~5). La section
suivante montre que les vari\'et\'es toriques lisses 
et les r\'eseaux
hermitiens sont des exemples de vari\'et\'es sur ${\mathbb
F}_1$.

Dans le paragraphe 6 nous d\'efinissons, dans certains cas, la
fonction z\^eta $\z_X (s)$ d'une vari\'et\'e $X$ sur ${\mathbb
F}_1$, \`a partir du nombre de points de $X$ dans les extensions
finies de ${\mathbb F}_1$. Cette fonction $\z_X (s)$ est un
polyn\^ome de la variable $s$, qui est bien celui pr\'evu par
Manin dans son article \cite{M} sur le sujet.

Le dernier paragraphe propose des sp\'eculations sur l'image de
la $K$-th\'eorie alg\'ebrique de ${\mathbb F}_1$ dans celle de
${\mathbb Z}$. Il est naturel de penser que ce groupe est
l'image de l'homomorphisme $J$ d'Adams, qui prend ses valeurs
dans les groupes d'homotopie stable des sph\`eres. Un r\'esultat
de Totaro sur les vari\'et\'es toriques \cite{To} se r\'ev\`ele
tout \`a fait coh\'erent avec une interpr\'etation de cette
image de $J$ en termes d'extensions de motifs de Tate mixtes sur
${\mathbb F}_1$.

Ce travail n'est bien s\^ur qu'une tentative. Il serait
int\'eressant de trouver d'autres exemples de vari\'et\'es sur
${\mathbb F}_1$ et d\'emontrer d'avantage de propri\'et\'es,
quitte \`a modifier les d\'efinitions. Par exemple, on aimerait
disposer de produits fibr\'es dans la cat\'egorie des sch\'emas
sur ${\mathbb F}_1$. On signalera au cours du texte quelques
variantes possibles des d\'efinitions pr\'esent\'ees.

Une version pr\'eliminaire de cet article est le court preprint
\cite{So}. J'ai b\'en\'efici\'e durant ce travail de tr\`es
nombreuses discussions, avec notamment J.-B.~Bost, M.~Brou\'e,
P.~Cartier, A.~Connes, I.~Gelfand, H.~Gillet, M.~Kapranov,
M.~Kontsevich, L.~Lafforgue, Y.~Manin, J-P. Serre
et B.~Totaro, que je tiens
\`a remercier.

\section{Un historique}\label{sec1}

\subsection{ \ } J.~Tits parle dans \cite{T} \S~13 du ``corps de
caract\'eristique un''. Si $G$ est un groupe de Chevalley,
simple et simplement connexe, le groupe des points de $G$ dans
${\mathbb F}_1$ n'est autre que le groupe de Weyl $W$ de $G$~:
\begin{equation}
\label{eq1}
W = G ({\mathbb F}_1) \, .
\end{equation}
Tits montre qu'on peut associer \`a $G$ une ``g\'eom\'etrie''
finie, dont $W$ est le groupe d'automorphismes et dont les
propri\'et\'es sont comparables \`a celles du groupe fini $G
({\mathbb F}_q)$ pour tout corps fini ${\mathbb F}_q$, $q > 1$.
Par exemple, l'hexagone est la g\'eom\'etrie associ\'ee \`a $G_2
({\mathbb F}_1)$.

La th\'eorie  des repr\'esentations fournit aussi de nombreux
exemples justifiant la formule (\ref{eq1}). R.~Steinberg avait
d\'ej\`a construit dans \cite{St} des repr\'esentations
irr\'eductibles de ${\rm GL}_N ({\mathbb F}_q)$ de fa\c con
parall\`ele \`a celle des repr\'esentations du groupe
sym\'etrique de $N$ lettres $\Si_N$. Leurs caract\`eres sont des
$q$-analogues des formules connues pour le groupe sym\'etrique.
Par exemple, si $N = N_1 + \cdots + N_r$, $N_i \geq 1$, est une
partition de l'entier $N$ et si $P \sbs {\rm SL}_N$ est le
sous-groupe parabolique standard associ\'e \`a cette partition,
le cardinal de l'ensemble fini ${\rm SL}_N ({\mathbb F}_q) /
P({\mathbb F}_q)$ est
$$
\# \, ({\rm SL}_N ({\mathbb F}_q) / P({\mathbb F}_q)) =\{ N \}
/ \{ N_1 \} \{ N_2 \} \ldots \{ N_r \} \, ,
$$
avec
$$
\{ n \}= \prod_{i=1}^{n} \ [i]
$$
et
$$
[n]= q^{n-1} + q^{n-2} + \cdots + 1 \, .
$$
Une telle formule, vraie si $q > 1$, demeure valable quand $q=1$
si l'on pose
$$
\Si_N = {\rm SL}_N ({\mathbb F}_1)
$$
(cf. (\ref{eq1})) et
$$
P ({\mathbb F}_1) = \prod_{i=1}^{r} \Si_{N_i} \, .
$$

Si $G$ est un groupe de Chevalley quelconque, 
les repr\'esentations
``unipotentes'' de $G ({\mathbb F}_q)$ sont des $q$-analogues
des repr\'esentations de $W$. Les travaux plus r\'ecents de
M.~Brou\'e, G.~Malle et J.~Michel \cite{B1}, \cite{B2},
 \cite{B3} montrent
qu'il est possible d'\'etendre cette analogie en faisant de $q$
une variable abstraite.

\subsection{ \ } Y.~Manin aborde d'un autre point de vue la
discussion du corps \`a un \'el\'ement ${\mathbb F}_1$
\cite{Ma}. L'analogie entre corps de nombres et corps de
fonctions $[W]$, dont est issue la th\'eorie d'Arakelov, 
fait du sch\'ema ${\rm Spec} ({\mathbb Z})$ l'analogue d'une
courbe affine sur un corps. On peut voir dans ${\mathbb F}_1$ le
corps de d\'efinition de la courbe ${\rm Spec} ({\mathbb Z})$.
Manin demande alors quel sens donner au produit fibr\'e
\begin{equation}
\label{eq2}
{\rm Spec} ({\mathbb Z}) \build\ts_{{\rm Spec} ({\mathbb F}_1)}
^{} {\rm Spec} ({\mathbb Z}) \, .
\end{equation}
Une telle question est naturelle si l'on se souvient du r\^ole
jou\'e par la surface $C \build\ts_{{\mathbb F}_q}^{} C$ dans la
preuve par Weil de l'hypoth\`ese de Riemann pour une courbe
projective et lisse $C$ sur le corps fini ${\mathbb F}_q$.

Manin propose aussi de d\'evelopper une g\'eom\'etrie
alg\'ebrique sur le corps ${\mathbb F}_1$ (\cite{Ma}, 1.7), ainsi
qu'une th\'eorie des motifs et des fonctions z\^etas
associ\'ees. Par exemple, l'espace projectif ${\mathbb
P}_{{\mathbb F}_1}^d$ de dimension $d$ sur ${\mathbb F}_1$
devrait avoir pour fonction z\^eta le polyn\^ome
\begin{equation}
\label{eq3}
\z_{{\mathbb P}_{{\mathbb F}_1}^d} = (s) s(s-1) (s-2) \ldots
(s-d) \, .
\end{equation}

\subsection{ \ } Ces sp\'eculations ont \'et\'e poursuivies par
Smirnov \cite{S} et Kapranov-Smirnov \cite{KS} \cite{K} (non
publi\'es). Ceux-ci proposent une alg\`ebre lin\'eaire sur
${\mathbb F}_1$ (un espace vectoriel lin\'eaire sur ${\mathbb
F}_1$ est un ensemble fini point\'e) et expliquent en ces termes
la loi de r\'eciprocit\'e de Gauss (voir \cite{A} pour le cas
g\'eom\'etrique). Ils proposent aussi que, si $n \geq 1$, le
corps ${\mathbb F}_1$ a une extension finie de degr\'e $n$, dont
les \'el\'ements inversibles sont les racines de l'unit\'e
d'ordre $n$.

\section{Pr\'eliminaires}\label{sec2}

\subsection{ \ } Notre objectif est de d\'efinir les
vari\'et\'es alg\'ebriques sur le corps \`a un \'el\'ement.

Notons d'abord que l'on ne s'attend pas \`a voir figurer ${\rm
Spec} ({\mathbb Z})$ parmi ces vari\'et\'es. En effet sa
fonction z\^eta, la fonction z\^eta de Riemann, a un nombre
infini de z\'eros, ce qui indique que ${\rm Spec} ({\mathbb Z})$
n'est pas de type fini sur ${\rm Spec} ({\mathbb F}_1)$. Et nous
ne chercherons pas \`a donner un sens au produit fibr\'e
(\ref{eq2}).

Par contre, si $X$ est une vari\'et\'e (de type fini) sur
${\mathbb F}_1$, on s'attend \`a ce que $X \build\ot_{{\mathbb
F}_1}^{} {\mathbb Z}$ soit une vari\'et\'e alg\'ebrique (de type
fini) sur ${\mathbb Z}$. Si de plus $X$ est lisse sur ${\mathbb
F}_1$, la vari\'et\'e $X \build\ot_{{\mathbb F}_1}^{} {\mathbb
Z}$ aura bonne r\'eduction partout, de r\'eduction en $p$ la
vari\'et\'e alg\'ebrique $X \build\ot_{{\mathbb F}_1}^{} {\mathbb
F}_p$, lisse que ${\mathbb F}_p$. Ceci nous conduit \`a la
question suivante~:

\medskip

\noindent {\bf Question 1.} {\it Quelles sont les vari\'et\'es
sur ${\mathbb Z}$ obtenues par extension des scalaires de
${\mathbb F}_1$ \`a ${\mathbb Z}$~?}

\medskip

On aimerait par exemple que les vari\'et\'es toriques ou les
vari\'et\'es de drapeaux d'un groupe de Chevalley puissent
\^etre d\'efinies sur ${\mathbb F}_1$.

\subsection{ \ } Un point de d\'epart pour nos d\'efinitions est
cette d\'efinition tr\`es \'economique des sch\'emas~: un
sch\'ema est un foncteur covariant de la cat\'egorie des anneaux
vers celle des ensembles qui est localement repr\'esentable par
un anneau. On trouvera dans \cite{DG} la preuve que cette
d\'efinition est \'equivalente \`a la d\'efinition usuelle. Ceci
sugg\`ere de d\'efinir les vari\'et\'es sur ${\mathbb F}_1$
comme foncteurs d'une cat\'egorie d'anneaux vers celle des
ensembles finis.

\subsection{ \ } Puisqu'\`a une vari\'et\'e $X$ sur ${\mathbb
F}_1$ on souhaite associer une vari\'et\'e alg\'ebrique $X
\build\ot_{{\mathbb F}_1}^{} {\mathbb Z}$ sur les entiers, il
convient de comprendre quel est le foncteur associ\'e \`a
l'extension des scalaires d'une vari\'et\'e alg\'ebrique.

Soit $k$ un corps, ${\mathcal A}_k$ la cat\'egorie des
$k$-alg\`ebres unitaires commutatives, et $\Om$ un objet de
${\mathcal A}_k$. Si $R$ est un objet de ${\mathcal A}_k$, on
note $R_{\Om} R \build\ot_{k}^{} \Om$ son extension des
scalaires de $k$ \`a $\Om$. De m\^eme, si $X$ est un sch\'ema
sur $k$, on pose $X_{\Om} X \build\ot_{k}^{} \Om$.

Notons ${\mathcal E}ns$ la cat\'egorie des ensembles. Si $X$ est
un sch\'ema sur $k$ on d\'esigne par
$$
\underline X : {\mathcal A}_k \ra {\mathcal E}ns
$$
le foncteur covariant qui \`a $R$ associe $X(R)$. De m\^eme, si
$S$ est un sch\'ema sur $\Om$, on d\'esigne par
$$
\underline S : {\mathcal A}_k \ra {\mathcal E}ns
$$
le foncteur covariant qui \`a $R$ associe $S (R_{\Om})$. La
proposition suivante montre que le foncteur $\underline X_{\Om}$
v\'erifie une propri\'et\'e universelle.

\medskip

\noindent {\bf Proposition 1.} i) {\it Pour tout objet $R$ de
${\mathcal A}_k$ et tout sch\'ema $X$ sur $k$ il existe une
inclusion canonique naturelle
$$
X(R) \sbs X_{\Om} (R_{\Om}) \, .
$$
On note $i : \underline X \ra \underline X_{\Om}$ la
transformation naturelle ainsi d\'efinie.}

\smallskip

\noindent ii) {\it Si $S$ est un sch\'ema sur $\Om$ et si
$$
\vp : \underline X \ra \underline S
$$
est une transformation naturelle, il existe un unique morphisme
alg\'ebrique sur $\Om$
$$
\vp_{\Om} : X_{\Om} \ra S
$$
tel que la transformation compos\'ee
$$
\begin{CD}
\underline X @>{i}>> \underline X_{\Om}
@>{\underline\vp_{\Om}}>> \underline S
\end{CD}
$$
de $i$ avec la transformation induite par $\vp_{\Om}$
co{\"\i}ncide avec $\vp$, i.e. $\vp =\underline\vp_{\Om} \circ
i$.}

\medskip

\noindent {\bf Preuve.} Puisque $k$ est contenu dans $\Om$ et que
$X_{\Om} (R_{\Om})$ co{\"\i}ncide avec l'ensemble des points de
$X$ dans la $k$-alg\`ebre $R_{\Om}$, l'\'enonc\'e i) est clair.
Pour montrer ii) supposons d'abord que $X$ est affine et soit $A=
\G (X ,{\mathcal O}_X)$ la $k$-alg\`ebre des fonctions globales sur $X$.
L'identit\'e de $A$ d\'efinit un point ${\rm id}_A \in X(A)$,
dont l'image par $\vp$ est un morphisme alg\'ebrique sur $\Om$~:
$$
\vp_{\Om} \in {\rm Hom}_{\Om} (X_{\Om} , S) = X_{\Om} (A_{\Om})
\, .
$$
Si $R$ est une $k$-alg\`ebre et $f \in X(R)$ un morphisme de $A$
vers $R$, on a, par fonctorialit\'e de $\vp$,
$$
\vp (f) =f^* (\vp ({\rm id}_A)) =f^* (\vp_{\Om})= \vp_{\Om}
\circ f_{\Om} \, .
$$
Cela montre que $\vp= \underline\vp_{\Om} \circ i$.

Dans le cas g\'en\'eral, soit $X \, =\build\cup_{j \in J}^{}
X_j$ un recouvrement ouvert 
de $X$ par des vari\'et\'es affines et
$$
\vp : \underline X \ra \underline S
$$
une transformation naturelle. La restriction de $\underline\vp$
\`a $\underline X_j$ (resp. $\underline{X_j \cap X_{j'}}$) est
induite par un morphisme alg\'ebrique $\vp_{j\Om} : X_{j\Om} \ra
S$ (resp. $\vp_{jj'\Om} : (X_j \cap X_{j'})_{\Om} \ra S$) et,
par fonctorialit\'e et unicit\'e, la restriction de $\vp_{j\Om}$
\`a $(X_j \cap X_{j'})_{\Om}$ co{\"\i}ncide avec $\vp_{jj'\Om}$
quels que soient les indices $j$ et $j'$ dans $J$. Par
cons\'equent, les morphismes $\vp_{j\Om}$, $j \in J$, se
recollent pour donner un morphisme
$$
\vp_{\Om} : X \ra S
$$
alg\'ebrique sur $\Om$. Si $R$ est une $k$-alg\`ebre et $f :
{\rm Spec} (R) \ra X$ un point de $X(R)$, consid\'erons les
$k$-sch\'emas affines $U_j =f^{-1} (X_j)$, $j \in J$. La
restriction de $(\vp_{\Om} \circ i)(f)$ \`a $U_{j\Om}$
co{\"\i}ncide par d\'efinition avec celle de $\vp_{j\Om} \circ
f$, c'est-\`a-dire avec la restriction de $\vp (f)$. Comme les
$U_{j\Om}$, $j \in J$, forment un recouvrement ouvert de
$U_{\Om}$, on en conclut que $\vp_{\Om} \circ i =\vp$. \hfill
q.e.d.

\subsection{ \ } On voudrait qu'un \'enonc\'e tel que la
Proposition~1 soit vrai quand $k = {\mathbb F}_1$ et $\Om =
{\mathbb Z}$. Mais cela suppose qu'on sache d'abord r\'epondre
\`a la question suivante~:

\medskip

\noindent {\bf Question 2.} {\it Quels sont les anneaux obtenus
par extension des scalaires de ${\mathbb F}_1$ \`a ${\mathbb
Z}$~?}

\medskip

Or, d'apr\`es Kapranov et Smirnov, pour tout entier $n \geq 1$,
le corps ${\mathbb F}_1$ poss\`ede une extension ${\mathbb
F}_{1^n}$ de degr\'e $n$, obtenue par l'adjonction des racines
$n$-i\`emes de l'unit\'e. La ${\mathbb Z}$-alg\`ebre ${\mathbb
F}_{1^n} \build\ot_{{\mathbb F}_1}^{} {\mathbb Z}$ est alors de
rang $n$ sur ${\mathbb Z}$. Cela conduit \`a poser
\begin{equation}
\label{eq3bis}
{\mathbb F}_{1^n} \build\ot_{{\mathbb F}_1}^{} {\mathbb Z} =
{\mathbb Z} [T] / (T^n - 1) \, .
\end{equation}
Cet anneau $R_n$ est donc une des r\'eponses \`a la Question~2,
et d'apr\`es la Proposition~1, si $X$ est une vari\'et\'e sur
${\mathbb F}_1$, il existe une inclusion naturelle
$$
X ({\mathbb F}_{1^n}) \sbs (X \build\ot_{{\mathbb F}_1}^{}
{\mathbb Z}) (R_n) \, .
$$

Plus g\'en\'eralement, nous admettrons que l'extension des
scalaires de ${\mathbb F}_1$ \`a ${\mathbb Z}$ induit une
\'equivalence de cat\'egorie entre une cat\'egorie de ${\mathbb
F}_1$-alg\`ebres et la cat\'egorie ${\mathcal R}$ des anneaux
finis et plats sur ${\mathbb Z}$, c'est-\`a-dire celle des
anneaux $R$ dont le groupe associ\'e est un ${\mathbb
Z}$-module libre de type fini. Une vari\'et\'e $X$ sur ${\mathbb
F}_1$ doit donc d\'efinir un foncteur covariant
$\underline X$ de la cat\'egorie
${\mathcal R}$ dans celle des ensembles, contenu dans le
foncteur qui \`a $R$ associe l'ensemble $(X \build\ot_{{\mathbb
F}_1}^{} {\mathbb Z}) (R)$.

Pour d\'efinir les vari\'et\'es sur ${\mathbb F}_1$ nous
proc\'ederons en deux temps. La cat\'egorie ${\mathcal A}$ des
vari\'et\'es affines sur ${\mathbb F}_1$ sera d\'efinie par une
propri\'et\'e universelle inspir\'ee de la Proposition~1 pour
des foncteurs de ${\mathcal R}$ vers ${\mathcal E}ns$. Celle des
vari\'et\'es sur ${\mathbb F}_1$ sera alors obtenue en
consid\'erant une propri\'et\'e universelle pour des foncteurs
de ${\mathcal A}$ vers  ${\mathcal E}ns$. 

\section{D\'efinitions} \label{sec3}

Rappelons que ${\mathcal R}$ est la cat\'egorie des anneaux
finis et plats sur ${\mathbb Z}$ et ${\mathcal E}ns$ celle des
ensembles.

\subsection{D\'efinition 1.} Un {\it truc sur} ${\mathbb F}_1$ est
le couple $X =(\underline X , \hbox{\LARGE $a$}_X)$ d'un foncteur
covariant
$$
\underline X : {\mathcal R} \ra {\mathcal E}ns
$$
et d'une ${\mathbb C}$-alg\`ebre $\hbox{\LARGE $a$}_X$. Pour tout
morphisme d'anneaux unitaires $\s : R \ra {\mathbb C}$, $R \in
{\rm Ob} ({\mathcal R})$, et pour tout \'el\'ement $x \in
\underline X (R)$ on suppose de plus donn\'e un morphisme
d'alg\`ebre (dit ``d'\'evaluation'')
$$
e_{x,\s} : \hbox{\LARGE $a$}_X \ra {\mathbb C} \, .
$$
Si $f : R' \ra R$ est un morphisme de ${\mathcal R}$ et si $y \in
\underline X (R')$, l'\'egalit\'e suivante doit \^etre
satisfaite~:
\begin{equation}
\label{eq4}
e_{f(y),\s}= e_{y , \s \circ f}
\end{equation}
Pour tout morphisme $\s : R \ra {\mathbb C}$.

\bigskip

\noindent {\bf Variantes.} i) On pourrait remplacer ${\mathcal R}$
par la sous-cat\'egorie pleine engendr\'ee par les anneaux $R_n$,
$n \geq 1$, et leurs produits tensoriels.

\noindent ii) L'id\'ee d'ajouter \`a $\underline X$ une ``donn\'ee
topologique \`a l'infini'' est due \`a J.-B.~Bost. Il n'y a pas
lieu de supposer que $\hbox{\LARGE $a$}_X$ est commutative. Un
objectif est de faire le lien avec la th\'eorie d'A.~Connes (voir
par exemple \cite{C}), un truc $X$ sur ${\mathbb F}_1$ disposant
d'une extension au point \`a l'infini de ${\rm Spec} ({\mathbb
Z})$, essentiellement donn\'ee par l'alg\`ebre $\hbox{\LARGE
$a$}_X$.
Par ignorance, nous resterons tr\`es \'evasifs sur les
propri\'et\'es requises sur cette ${\mathbb C}$-alg\`ebre. 

\subsection{D\'efinition 2.} i) Un truc $X =(\underline X ,
\hbox{\LARGE $a$}_X)$ sur ${\mathbb F}_1$ est {\it fini} quand
tous les ensembles $\underline X (R)$, $R \in {\rm Ob} ({\mathcal
R})$, sont finis.

\smallskip

\noindent ii) Un {\it morphisme} $\vp : X \ra Y$ entre deux trucs
sur ${\mathbb F}_1$ est la donn\'ee d'une transformation naturelle
$$
\underline\vp : \underline X \ra \underline Y
$$
et d'un morphisme d'alg\`ebres
$$
\vp^* : \hbox{\LARGE $a$}_Y \ra \hbox{\LARGE $a$}_X
$$
tels que, si $R \in {\rm Ob} ({\mathcal R})$, si $\s : R \ra
{\mathbb C}$ est un morphisme d'anneaux unitaires, et si $x \in
\underline X (R)$, on a
\begin{equation}
\label{eq5}
e_{\s , \underline\vp (x)} (\a)= e_{\s , x} (\vp^* (\a))
\end{equation}
pour toute fonction $\a \in \hbox{\LARGE $a$}_Y$.

\smallskip

\noindent iii) Un morphisme
$$
\vp =(\underline\vp , \vp^*) : X \ra Y
$$
est une {\it immersion} si $\vp^*$ est injectif et si
l'application
$$
\underline\vp : \underline X (R) \ra \underline Y (R)
$$
est injective quel que soit l'anneau $R \in {\rm Ob} ({\mathcal
R})$.

\bigskip
\bigskip

Si $\vp =(\underline\vp , \vp^*) : X \ra Y$ et $\psi =
(\underline \psi , \psi^*) : Y \ra Z$ sont deux morphismes, leur
compos\'e est le couple
$$
\psi \circ \vp =(\underline\psi \circ \underline\vp , \vp^*
\circ \psi^*) : X \ra Z \, .
$$
On note ${\mathcal T}$ la cat\'egorie des trucs sur ${\mathbb
F}_1$.

\subsection{ \ } Soit ${\mathcal V}_{\mathbb Z}$ la cat\'egorie
des vari\'et\'es sur ${\mathbb Z}$, c'est-\`a-dire les sch\'emas
de type fini sur ${\rm Spec} ({\mathbb Z})$. Si $V$ et $W$ sont
deux objets de ${\mathcal V}_{\mathbb Z}$, on d\'esigne par ${\rm
Hom}_{\mathbb Z} (V,W)$ l'ensemble des morphismes alg\'ebriques de
$V$ vers $W$.

Si $V$ est un objet de ${\mathcal V}_{\mathbb Z}$, on lui associe
un truc $V =(\underline V , {\mathcal O} (V_{\mathbb C}))$ sur
${\mathbb F}_1$ de la fa\c con suivante. Si $R \in {\rm Ob}
({\mathcal R})$ on pose
$$
\underline V (R)= {\rm Hom}_{\mathbb Z} ({\rm Spec} (R) , V)
$$
et si $f : R' \ra R$ est une fl\`eche de ${\mathcal R}$,
l'application
$$
\underline V (f) : \underline V (R') \ra \underline V (R)
$$
est la composition avec $f$. La ${\mathbb C}$-alg\`ebre ${\mathcal
O} (V_{\mathbb C})$ est celle des fonctions globales (i.e. les
sections du faisceau canonique) sur la vari\'et\'e $V_{\mathbb C}=
V \build\ot_{\mathbb Z}^{} {\mathbb C}$. Si $\s : R \ra
{\mathbb C}$ est un morphisme d'anneaux unitaires et si $x \in
\underline V (R)$, l'image de $x$ par $\s$ est un point complexe
$\s (x)$ de $V_{\mathbb C}$. On note
$$
e_{x,\s} :{\mathcal O} (V_{\mathbb C}) \ra {\mathbb C}
$$
l'\'evaluation au point $\s (x)$ des fonctions sur $V_{\mathbb
C}$. La formule (\ref{eq4}) est \'evidemment v\'erifi\'ee. Tout
morphisme alg\'ebrique $f : V \ra W$ induit un morphisme dans
${\mathcal T}$, \'egalement not\'e $f$.

\subsection{D\'efinition 3.}

Une {\it vari\'et\'e affine sur} ${\mathbb F}_1$ est un truc fini
$X$ sur ${\mathbb F}_1$ tel qu'il existe une vari\'et\'e
alg\'ebrique affine $X_{\mathbb Z}$ sur ${\mathbb Z}$ et une
immersion $i : X \ra X_{\mathbb Z}$ dans ${\mathcal T}$
v\'erifiant la propri\'et\'e suivante.

Quels que soient la vari\'et\'e affine $V \in {\rm Ob} ({\mathcal
V}_{\mathbb Z})$ et le morphisme de ${\mathcal T}$
$$
\vp : X \ra V \, ,
$$
il existe un unique morphisme alg\'ebrique
$$
\vp_{\mathbb Z} : X_{\mathbb Z} \ra V
$$
tel que $\vp = \vp_{\mathbb Z} \circ i$.

\bigskip
\bigskip

On voit que la vari\'et\'e $X_{\mathbb Z}$ dans ${\mathcal
V}_{\mathbb Z}$ est uniquement d\'etermin\'ee par $X$. On la note
aussi $X \build\ot_{{\mathbb F}_1}^{} {\mathbb Z}$. Par
ailleurs, si $X$ et $Y$ sont deux vari\'et\'es affines sur
${\mathbb F}_1$ et si $f : X \ra Y$ est un morphisme de
${\mathcal
T}$, il existe un unique morphisme $f_{\mathbb Z} \in {\rm
Hom}_{\mathbb Z} (X_{\mathbb Z} , Y_{\mathbb Z})$ qui rend
commutatif le diagramme
$$
\begin{CD}
X_{\mathbb Z} @>{f_{\mathbb Z}}>> Y_{\mathbb Z} \\
@AAA @AAA \\
X @>{f}>> Y
\end{CD}
$$
On notera ${\mathcal A}$ la sous-cat\'egorie pleine de ${\mathcal
T}$ dont les objets sont les vari\'et\'es affines.

\subsection{D\'efinition 4.}

\noindent i) Un {\it objet sur} ${\mathbb F}_1$ est la donn\'ee $X=
(\underline{\underline X} , \hbox{\LARGE $a$}_X)$ d'un foncteur
contravariant
$$
\underline{\underline X} : {\mathcal A} \ra {\mathcal E}ns
$$
et d'une ${\mathbb C}$-alg\`ebre $\hbox{\LARGE $a$}_X$ ainsi que
d'un morphisme d'alg\`ebres
$$
e_x : \hbox{\LARGE $a$}_X \ra \hbox{\LARGE $a$}_A
$$
pour objet $A$ de ${\mathcal A}$ et tout \'el\'ement $x$ de
$\underline{\underline X} (A)$. Si $f : A \ra B$ est un morphisme
de ${\mathcal A}$ et $x \in \underline{\underline X} (B)$, on
suppose de plus que l'\'egalit\'e
\begin{equation}
\label{eq6}
e_{f^* (x)} = f^* \circ e_x
\end{equation}
est v\'erifi\'ee, avec $f^* (x)= \underline{\underline X} (f)(x)
\in \underline{\underline X} (A)$.

\smallskip

\noindent ii) Un objet $X$ sur ${\mathbb F}_1$ est {\it fini} si
$\underline{\underline X} ({\rm Spec} \, R)$ est fini quel que soit
$R \in {\rm Ob} ({\mathcal R})$ (d'apr\`es la
 Proposition~2 ii) ci-dessous, la
cat\'egorie ${\mathcal R}^{\rm opp}$ est contenue dans ${\mathcal
A}$).

\smallskip

\noindent iii) Un {\it morphisme} $\vp : X \ra Y$ entre objets sur
${\mathbb F}_1$ est la donn\'ee d'une transformation naturelle
$$
\underline{\underline \vp} : \underline{\underline X} \ra
\underline{\underline Y}
$$
et d'un morphisme d'alg\`ebres $\vp^* : \hbox{\LARGE $a$}_Y \ra
\hbox{\LARGE $a$}_X$ tels que, si $A \in {\rm Ob} ({\mathcal A})$
et $x \in \underline{\underline X} (A)$, on ait
\begin{equation}
\label{eq7}
e_{\underline{\underline \vp} (x)}=e_x \circ \vp^* \, .
\end{equation}

\noindent iv) On dit que le morphisme $\vp$ est une {\it
immersion} si $\underline{\underline \vp}$ et $\vp^*$ sont
injectifs.

\bigskip
\bigskip

Le compos\'e de deux morphismes est d\'efini de la fa\c con
\'evidente, et l'on note ${\mathcal O}$ la cat\'egorie des objets
sur ${\mathbb F}_1$.

\subsection{ \ } Si $V \in {\rm Ob} ({\mathcal V}_{\mathbb Z})$,
on lui associe comme suit un objet $V =(\underline{\underline V}
, {\mathcal O} (V_{\mathbb C}))$ sur ${\mathbb F}_1$.

Si $A \in {\rm Ob} ({\mathcal A})$ on pose
$$
\underline{\underline V} (A) ={\rm Hom}_{\mathbb Z} (A_{\mathbb
Z} , V)
$$
et si $f : A \ra B$ est une fl\`eche de ${\mathcal A}$ on
d\'esigne par $\underline{\underline V} (f)$ la composition avec
$f_{\mathbb Z}$. Si $x \in \underline{\underline V} (A)$,
l'\'evaluation $e_x$ est le morphisme compos\'e
$$
\begin{CD}
{\mathcal O} (V_{\mathbb C}) @>{x^*}>> {\mathcal O} (A_{\mathbb
C}) @>{i^*}>>  \hbox{\LARGE $a$}_A \, ,
\end{CD}
$$
o\`u $i : A \ra A_{\mathbb Z}$ est l'inclusion associ\'ee \`a $A$.

En associant \`a $f \in {\rm Hom}_{\mathbb Z} (V,W)$ le morphisme
de composition avec $f$ et celui d'image inverse $f^* : {\mathcal
O} (W_{\mathbb C}) \ra {\mathcal O} (V_{\mathbb C})$, on obtient
ainsi un foncteur
$$
{\mathcal V}_{\mathbb Z} \ra {\mathcal O} \, .
$$

\subsection{D\'efinition 5}

Une {\it vari\'et\'e sur} ${\mathbb F}_1$ est la donn\'ee d'un
objet $X$ sur ${\mathbb F}_1$ tel qu'il existe une vari\'et\'e
alg\'ebrique $X_{\mathbb Z}$ sur ${\mathbb Z}$ et une immersion $i
: X \ra X_{\mathbb Z}$ de ${\mathcal O}$ ayant la propri\'et\'e
suivante. Pour toute vari\'et\'e $V \in {\rm Ob} ({\mathcal
V}_{\mathbb Z})$ et tout morphisme
$$
\vp : X \ra V
$$
dans ${\mathcal O}$, il existe un unique morphisme alg\'ebrique
$$
\vp_{\mathbb Z} : X_{\mathbb Z} \ra V
$$
tel que $\vp =\vp_{\mathbb Z} \circ i$.

\bigskip
\bigskip

La vari\'et\'e $X_{\mathbb Z} \in {\rm Ob} ({\mathcal V}_{\mathbb
Z})$ est uniquement d\'etermin\'ee (\`a isomorphisme unique pr\`es) par
la vari\'et\'e $X$ sur ${\mathbb F}_1$. On la note aussi $X
\build\ot_{{\mathbb F}_1}^{} {\mathbb Z}$ et on l'appelle {\it
extension des scalaires de} $X$ de ${\mathbb F}_1$ \`a ${\mathbb
Z}$. Si $X$ et $Y$ sont deux vari\'et\'es sur ${\mathbb F}_1$ et
$f : X \ra Y$ est un morphisme de ${\mathcal O}$, il existe un
unique morphisme alg\'ebrique $f_{\mathbb Z} \in {\rm
Hom}_{\mathbb Z} (X_{\mathbb Z}, Y_{\mathbb Z})$ qui induit $f$
sur $X$ (cf.~3.4). Autrement dit, si ${\mathcal V}$ est la
sous-cat\'egorie pleine de ${\mathcal O}$ dont les objets sont les
vari\'et\'es sur ${\mathbb F}_1$, l'extension des scalaires de
${\mathbb F}_1$ \`a ${\mathbb Z}$ est un foncteur fid\`ele de
${\mathcal V}$ vers ${\mathcal V}_{\mathbb Z}$.

\subsection{Variantes}

\subsubsection{ \ } Pour \'eviter le choix trop trivial
$\hbox{\LARGE $a$}_X ={\mathcal O} (X_{\mathbb C})$ (o\`u
$X_{\mathbb C} =X_{\mathbb Z} \build\ot_{\mathbb Z}^{}
{\mathbb C})$ dans la d\'efinition d'une vari\'et\'e $X$ sur
${\mathbb F}_1$ (cf. Proposition~4 ci-dessous) on pourrait par
exemple imposer (dans les d\'efinitions 3 et 5) que l'alg\`ebre
$\hbox{\LARGE $a$}_X$ est une alg\`ebre de Banach commutative. Ce
sera le cas pour les exemples discut\'es dans 
la section~5
(\cite{R}, 18.11).

\subsubsection{ \ } Comme l'a sugg\'er\'e R.~Pink, pour tout
entier $n \geq 1$, on peut aussi d\'efinir les vari\'et\'es sur
${\mathbb F}_{1^n}$ en rempla\c cant dans les d\'efinitions
pr\'ec\'edentes la cat\'egorie ${\mathcal R}$ (resp. ${\mathcal
V}_{\mathbb Z}$) par la sous-cat\'egorie des 
$R_n$-alg\`ebres finies et plates (resp. par la cat\'egorie des
sch\'emas de type fini sur ${\rm Spec} (R_n)$).

\section{Quelques propri\'et\'es}\label{sec4}

\subsection{ \ } La cat\'egorie oppos\'ee \`a ${\mathcal R}$ est
contenue dans ${\mathcal A}$.

\medskip

\noindent {\bf Proposition 2.} i) {\it Si $R \in {\rm Ob}
({\mathcal R})$ et $X \in {\rm Ob} ({\mathcal T})$, l'ensemble
$\underline X (R)$ est canoniquement isomorphe \`a ${\rm
Hom}_{\mathcal T} ({\rm Spec} (R) , X)$.}

\smallskip

\noindent ii) {\it Si $R \in {\rm Ob} ({\mathcal R})$, le truc sur
${\mathbb F}_1$ associ\'e \`a ${\rm Spec} (R)$ est une vari\'et\'e
affine sur ${\mathbb F}_1$ dont l'extension \`a ${\mathbb Z}$
co{\"\i}ncide avec ${\rm Spec} (R)$. On obtient ainsi un foncteur
contravariant pleinement fid\`ele de ${\mathcal R}$ dans
${\mathcal A}$.}

\medskip

\noindent {\bf Preuve.} i) Soient $R \in {\rm Ob} ({\mathcal R})$
un anneau fini et plat sur ${\mathbb Z}$ et ${\rm Spec} (R)$ le
sch\'ema associ\'e. Si $X \in {\rm Ob} ({\mathcal T})$ et si
$$
\vp : {\rm Spec} (R) \ra X
$$
est un morphisme de ${\mathcal T}$, l'image par $\underline\vp$ de
l'application identique
$$
{\rm id}_R \in {\rm Hom}_{\mathbb Z} ({\rm Spec} (R), {\rm Spec}
(R))
$$
est un \'el\'ement $\underline\vp ({\rm id}_R) \in \underline X
(R)$.

Inversement, \'etant donn\'e $x \in \underline X (R)$, pour tout
morphisme de ${\mathcal R}$
$$
f \in {\rm Hom}_{\mathcal R} (R,R') = {\rm Hom}_{\mathbb Z} ({\rm
Spec} (R'), {\rm Spec} (R)) \, ,
$$
l'image de $x$ par l'application
$$
\underline X (f) : \underline X (R) \ra \underline X (R')
$$
d\'efinit un \'el\'ement $f(x) \in \underline X (R')$ et l'on
obtient ainsi une transformation naturelle
$$
\underline x : \underline{\rm Spec}(R) \ra \underline X \, .
$$
Par ailleurs, si $\Si$ est l'ensemble (fini) des morphismes
unitaires $\s : R \ra {\mathbb C}$, on dispose d'isomorphismes
canoniques
$$
{\mathcal O} ({\rm Spec} (R)_{\mathbb C}) = R
\build\ot_{\mathbb Z}^{} {\mathbb C} = {\mathbb C}^{\Si} \, .
$$
La collection des morphismes d'\'evaluation en $x$
$$
e_{x , \s} : \hbox{\LARGE $a$}_X \ra {\mathbb C} \, , \qquad \s
\in \Si \, ,
$$
d\'efinit donc un morphisme d'alg\`ebres
$$
x^* : \hbox{\LARGE $a$}_X \ra R \build\ot_{\mathbb Z}^{} {\mathbb
C}
$$
qui v\'erifie la condition (\ref{eq5}) avec $\underline x$. On
obtient ainsi un morphisme $(\underline x , x^*)$ de ${\rm Spec}
(R)$ vers $X$ dans ${\mathcal T}$.

On v\'erifie que les applications $\vp \mpo \underline\vp ({\rm
id}_R)$ et $x \mpo (\underline x , x^*)$ sont des bijections
inverses entre ${\rm Hom}_{\mathcal T} ({\rm Spec} (R) , X)$ et
$\underline X (R)$.

Par exemple on a la formule
\begin{equation}
\label{eq8}
X ({\mathbb F}_{1^n}) = \underline X (R_n) \, ,
\end{equation}
o\`u le terme de gauche d\'esigne les morphismes dans ${\mathcal
T}$ de ${\rm Spec} ({\mathbb F}_{1^n})$ vers $X$.

\smallskip

\noindent ii) Si $R$ et $R'$ sont des anneaux de ${\mathcal R}$,
l'ensemble $\underline{{\rm Spec} (R)} (R')$ des morphismes
d'anneaux de $R$ vers $R'$ est fini. Par ailleurs, si $V$ est un
objet de ${\mathcal V}_{\mathbb Z}$ et $R$ un objet de ${\mathcal
R}$ on sait, d'apr\`es i), que
$$
\underline V (R) ={\rm Hom}_{\mathcal T} ({\rm Spec} (R) , V) \,
.
$$
Comme $\underline V (R)= {\rm Hom}_{\mathbb Z} ({\rm Spec} (R) ,
V)$ par d\'efinition de $V$, on voit que tout morphisme de ${\rm
Spec} (R)$ vers $V$ dans ${\mathcal T}$ est alg\'ebrique. Cela
montre que ${\rm Spec} (R)$ est une vari\'et\'e affine sur
${\mathbb F}_1$ dont l'extension \`a ${\mathbb Z}$ est ${\rm Spec}
(R)$.

En choisissant $V ={\rm Spec} (R')$, $R' \in {\rm Ob} ({\mathcal
R})$, dans 
l'argument pr\'ec\'edent on voit aussi que tout morphisme
$$
{\rm Spec} (R) \ra {\rm Spec} (R')
$$
dans ${\mathcal T}$ provient d'un unique morphisme $R' \ra R$ dans
${\mathcal R}$. Par cons\'equent le foncteur $R \mpo {\rm Spec}
(R)$ de ${\mathcal R}$ dans ${\mathcal A}$ est pleinement
fid\`ele.

\subsection{ \ } On d\'efinit un foncteur
$$
\ve : {\mathcal T} \ra {\mathcal O}
$$
en associant \`a un truc $X =(\underline X , \hbox{\LARGE
$a$}_X)$ le couple $(\underline{\underline X} , \hbox{\LARGE
$a$}_X)$ o\`u $\underline{\underline X}$ est le foncteur sur
${\mathcal A}$ repr\'esent\'e par $X$:
$$
\underline{\underline X} (A)= {\rm Hom}_{\mathcal T} (A,X) \, .
$$
Si $u \in \underline{\underline X} (A)$ l'\'evaluation en $u$ est
l'image inverse $u^*$:
$$
e_u = u^* : \hbox{\LARGE $a$}_X \ra \hbox{\LARGE $a$}_A \, .
$$
Ce foncteur $\ve$ va nous permettre de consid\'erer les
vari\'et\'es affines sur ${\mathbb F}_1$ comme des vari\'et\'es
sur ${\mathbb F}_1$~:

\medskip

\noindent {\bf Proposition 3.} i) {\it Le foncteur $\ve :
{\mathcal T} \ra {\mathcal O}$ est pleinement fid\`ele.}

\smallskip

ii) {\it L'image essentielle de ${\mathcal A}$ par $\ve$ est la
cat\'egorie des vari\'et\'es sur ${\mathbb F}_1$ dont l'extension
des scalaires \`a ${\mathbb Z}$ est affine.}

\medskip

\noindent {\bf Preuve.} i) Consid\'erons le foncteur
$$
\rho : {\mathcal O} \ra {\mathcal T}
$$
qui \`a l'objet $(\underline{\underline X} , \hbox{\LARGE $a$}_X)$
sur ${\mathbb F}_1$ associe le truc $(\underline X , \hbox{\LARGE
$a$}_X)$, o\`u $\underline X$ est la restriction \`a ${\mathcal
R}$ du foncteur $\underline{\underline X}$ (cf. Proposition~2 ii)).

Si $X$ est dans ${\mathcal T}$, le truc $\rho \circ \ve (X)$
v\'erifie
$$
\underline{\rho \circ \ve (X)} (R) ={\rm Hom}_{\mathcal T} ({\rm
Spec} (R) , X) = \underline X (R)
$$
(Proposition~2 i)). Par cons\'equent
\begin{equation}
\label{eq9}
\rho \circ \ve = {\rm id}_{\mathcal T} \, .
\end{equation}

De plus $\rho$ est l'adjoint \`a gauche
de $\ve$: si $Y$ est un truc sur
${\mathbb F}_1$ et $X$ un objet sur ${\mathbb F}_1$ il existe un
isomorphisme canonique et naturel
\begin{equation}
\label{eq10}
{\rm Hom}_{\mathcal T} (\rho (X),Y) = {\rm Hom}_{\mathcal O} (X ,
\ve (Y)) \, .
\end{equation}
En effet, si
$$
\psi : X \ra \ve (Y)
$$
est un morphisme de ${\mathcal O}$, le morphisme
$$
\rho (\psi) : \rho (X) \ra \rho \circ \ve (Y) =Y
$$
est un morphisme de ${\mathcal T}$.

Inversement, \'etant donn\'e un morphisme
$$
\vp : \rho (X) \ra Y
$$
de ${\mathcal T}$, pour tout objet $A$ de ${\mathcal A}$ et tout
$x \in \underline{\underline X} (A)$, si $R \in {\rm Ob}
({\mathcal R})$ et
$$
f \in \underline A (R)= {\rm Hom}_{\mathcal T} ({\rm Spec} (R),
A)
$$
(Prop.~2 i)), puisque $\underline{\underline X}$ est un foncteur
contravariant, on obtient un \'el\'ement $\underline{\underline X}
(f)(x) \in \underline{\underline X} (R)= \underline{\rho
(X)}(R)$, et donc un \'el\'ement
$$
\underline\vp (\underline{\underline X} (f)(x)) \in \underline Y
(R) \, .
$$
L'application $f \mpo \underline\vp (\underline{\underline X}
(f)(x))$ d\'efinit une transformation naturelle de $\underline A$
vers $\underline Y$ qui, jointe au morphisme d'alg\`ebres
$$
x^* \circ \vp^* : \hbox{\LARGE $a$}_Y \ra \hbox{\LARGE $a$}_A \, ,
$$
d\'efinit un \'el\'ement de
$$
{\rm Hom}_{\mathcal T} (A,Y)= \underline{\underline{\ve (Y)}} (A)
$$
qui d\'epend fonctoriellement de $x \in \underline{\underline X}
(A)$. D'o\`u une transformation naturelle
$$
\underline{\underline X} \ra \underline{\underline{\ve (Y)}}
$$
qui, jointe au morphisme d'alg\`ebres
$$
\vp^* :  \hbox{\LARGE $a$}_Y \ra \hbox{\LARGE $a$}_X \, ,
$$
fournit un morphisme $\a (\vp)$ de $X$ vers $\ve (Y)$. On
v\'erifie que les applications $\vp \mpo \a (\vp)$ et $\psi \mpo
\rho (\psi)$ sont des bijections inverses naturelles entre ${\rm
Hom}_{\mathcal T} (\rho (X) , Y)$ et ${\rm Hom}_{\mathcal O} (X ,
\ve (Y))$. Cela d\'emontre (\ref{eq10}).

\medskip

Il r\'esulte de (\ref{eq9}) et (\ref{eq10}) que $\ve : {\mathcal
T} \ra {\mathcal O}$ est pleinement fid\`ele.

\smallskip

\noindent ii) Consid\'erons une vari\'et\'e affine $X$ sur
${\mathbb F}_1$, $V \in {\rm Ob} ({\mathcal V}_{\mathbb Z})$ et
$$
\vp : \ve (X) \ra V
$$
un morphisme de ${\mathcal O}$. L'image par $\underline{\underline
\vp}$ du morphisme identique
$$
{\rm id}_X \in {\rm Hom}_{\mathcal T} (X,X) =
\underline{\underline{\ve (X)}} (X)
$$
est un morphisme alg\'ebrique $\vp_{\mathbb Z} \in {\rm
Hom}_{\mathbb Z} (X_{\mathbb Z} , V)$. Si $A$ est une vari\'et\'e
affine sur ${\mathbb F}_1$ et
$$
f \in {\rm Hom}_{\mathcal T} (A,X)= \underline{\underline X} (A)
\, ,
$$
on sait que $f$ est  par un morphisme alg\'ebrique
$f_{\mathbb Z} \in {\rm Hom}_{\mathbb Z} (A_{\mathbb Z} ,
X_{\mathbb Z})$ (D\'efinition~3) et donc
$$
\underline{\underline \vp} (f)= \underline{\underline \vp} (f^*
({\rm id}_X))= f_{\mathbb Z}^* (\underline{\underline \vp} ({\rm
id}_X))= f_{\mathbb Z}^* (\vp_{\mathbb Z}) =\vp_{\mathbb Z}
\circ f_{\mathbb Z}
$$
dans ${\rm Hom}_{\mathbb Z} (A_{\mathbb Z} , V)$. Il en r\'esulte
que $\vp$ est le morphisme induit par $\vp_{\mathbb Z}$ sur $\ve
(X)$, et donc $\ve (X)$ est une vari\'et\'e sur ${\mathbb F}_1$
dont $X_{\mathbb Z}$ est l'extension des scalaires \`a ${\mathbb
Z}$.

Inversement, supposons que $X$ soit une vari\'et\'e sur ${\mathbb
F}_1$ telle que $X_{\mathbb Z}$ soit affine. Soient $V \in {\rm
Ob} ({\mathcal V}_{\mathbb Z})$ une vari\'et\'e alg\'ebrique
affine et
$$
\vp : \rho (X) \ra V
$$
un morphisme de ${\mathcal T}$. Puisque $V$ est affine on a
$$
\ve (\underline V , {\mathcal O} (V_{\mathbb C})) =
(\underline{\underline V} , {\mathcal O} (V_{\mathbb C})) \, .
$$
Donc $\vp$ induit un morphisme
$$
\ve (\vp) : \ve \circ \rho (X) \ra V
$$
dans ${\mathcal O}$. L'isomorphisme d'adjonction
$$
{\rm Hom}_{\mathcal T} (\rho (X) , \rho (X)) ={\rm Hom}_{\mathcal O}
(X , \ve \circ \rho (X))
$$
associe \`a l'identit\'e de $\rho (X)$ un morphisme canonique
$$
X \ra \ve \circ \rho (X) \, .
$$
Puisque $X$ est une vari\'et\'e sur ${\mathbb F}_1$, le compos\'e
du morphisme pr\'ec\'edent avec $\ve (\vp)$ dans
${\mathcal O}$
est induit  par un morphisme alg\'ebrique
$$
\vp_{\mathbb Z} : X_{\mathbb Z} \ra V \, .
$$
En appliquant le foncteur $\rho$, on voit que $\vp_{\mathbb Z}
\circ i$ co{\"\i}ncide avec le morphisme
$$
\vp : \rho (X) \ra V
$$
dans ${\mathcal T}$. Par cons\'equent $\rho (X)$ est une
vari\'et\'e affine dont $X_{\mathbb Z}$ est l'extension des
scalaires \`a ${\mathbb Z}$. Cela d\'emontre ii).

\subsection{ \ } En g\'en\'eral une vari\'et\'e sur ${\mathbb Z}$
peut \^etre l'extension des scalaires \`a ${\mathbb Z}$ de
plusieurs vari\'et\'es sur ${\mathbb F}_1$.

\medskip

\noindent {\bf Proposition 4.} {\it Soit $X$ une vari\'et\'e sur
${\mathbb F}_1$. Supposons que l'immersion $i : X \ra X_{\mathbb
Z}$ soit la compos\'ee dans ${\mathcal O}$ d'une immersion $u : X
\ra Y$ et d'une immersion $j : Y \ra X_{\mathbb Z}$.

Supposons de plus que $X$ est affine ou que $\underline{\underline
u}$ est une \'equivalence. Alors $Y$ est une vari\'et\'e sur
${\mathbb F}_1$ telle que
$$
Y \build\otimes_{{\mathbb F}_1}^{} {\mathbb Z}= X_{\mathbb Z} \, .
$$
}

\medskip

\noindent {\bf Preuve.} Soient $V \in {\rm Ob} ({\mathcal
V}_{\mathbb Z})$ et
$$
\vp : Y \ra V
$$
un morphisme de ${\mathcal O}$. Puisque $X$ est une vari\'et\'e sur
${\mathbb F}_1$, la restriction de $\vp$ \`a $X$ est induite par un
morphisme alg\'ebrique $\vp_{\mathbb Z} \in {\rm Hom}_{\mathbb Z}
(X_{\mathbb Z} , V)$: $\vp_{\mathbb Z} \circ i =\vp \circ u$. Il
s'agit de v\'erifier que $\vp_{\mathbb Z} \circ j =\vp$.

Or le compos\'e des morphismes d'alg\`ebres
$$
\begin{CD}
{\mathcal O} (V_{\mathbb C}) @>{\vp_{\mathbb Z}^*}>> {\mathcal O}
(X_{\mathbb C}) @>{i^*}>> \hbox{\LARGE $a$}_X
\end{CD}
$$
co{\"\i}ncide avec
$$
\begin{CD}
{\mathcal O} (V_{\mathbb C}) @>{\vp^*}>> \hbox{\LARGE $a$}_Y
@>{u^*}>> \hbox{\LARGE $a$}_X \, .
\end{CD}
$$
Comme $i^* =u^* \circ j^*$ et comme $u^*$ est injectif, on en
d\'eduit que
\begin{equation}
\label{eq11}
\vp^*= j^* \circ \vp_{\mathbb Z}^* \, .
\end{equation}
Si le foncteur $\underline{\underline u} : \underline{\underline X}
\ra \underline{\underline Y}$ est une \'equivalence on a
$\underline{\underline \vp}_{\mathbb Z} \circ\underline{\underline
j}= \underline{\underline \vp}$ donc (\ref{eq11}) suffit \`a
montrer que $\vp_{\mathbb Z} \circ j= \vp$.

Dans le cas o\`u $X$ est affine, on peut supposer que $V$ est
affine (Proposition~3 ii)). Si $R \in {\rm Ob} ({\mathcal R})$, si
$x \in \underline{\underline Y} (R)$,
si $\a \in {\mathcal O} (V_{\mathbb C})$ et si $\s : R \ra {\mathbb
C}$ est un morphisme d'anneaux unitaires, on a, d'apr\`es
(\ref{eq7}) et (\ref{eq11})
\begin{equation}
\label{eq12}
e_{\s , \underline{\underline \vp} (x)} (\a) =e_{\s , x} (\vp^*
(\a))= e_{\s , x} (j^* \circ \vp_{\mathbb Z}^* (\a))= e_{\s ,
\underline{\underline \vp}_{\mathbb Z} 
\circ \underline{\underline j}(x)} (\a) \, .
\end{equation}
Comme $R$ est plat sur ${\mathbb Z}$, le morphisme canonique
$$
R \ra R \build\otimes_{\mathbb Z}^{} {\mathbb C} =R^{\Si}
$$
est injectif. Et comme $V$ est
affine les fonctions de ${\mathcal O} (V_{\mathbb C})$ s\'eparent
les points de $V({\mathbb C})$. Par cons\'equent les \'egalit\'es
(\ref{eq12}) avec $\s \in \Si$ et
 $\a \in {\mathcal O} (V_{\mathbb C})$
montrent que
$$
\underline{\underline \vp} (x) 
=\underline{\underline \vp}_{\mathbb Z}  \circ
\underline{\underline j} (x) \, .
$$
Par cons\'equent $\vp =\vp_{\mathbb Z} \circ j$. \hfill q.e.d.

\bigskip

\subsection{ \ } L'\'enonc\'e suivant permet de d\'efinir des
vari\'et\'es sur ${\mathbb F}_1$ par recollement.

\medskip

\noindent {\bf Proposition 5.} {\it Soient $V \in {\rm Ob}
({\mathcal V}_{\mathbb Z})$ et $V =\, \build\cup_{i \in I}^{} U_i$
un recouvrement ouvert fini de $V$. On suppose donn\'ees des
vari\'et\'es $X_i$, $i \in I$, et $X_{ij}$, $i \ne j$, sur
${\mathbb F}_1$ et des immersions
$$
X_{ij} \ra X_i \quad \hbox{et} \quad X_i \ra V
$$
dans ${\mathcal O}$ dont les extensions \`a 
${\mathbb Z}$ sont les
inclusions
$$
U_i \cap U_j \ra U_i \quad \hbox{et} \quad U_i \ra V \, .
$$
On suppose de plus que l'immersion compos\'ee
$$
X_{ij} \ra X_i \ra V
$$
co{\"\i}ncide avec
$$
X_{ij} \ra X_j \ra V
$$
si $i \ne j$. Pour tout objet $A$ de ${\mathcal A}$ on pose
$$
\underline{\underline X} (A) = \,
 \build\cup_{i}^{} \underline{\underline
X_i} (A)
$$
(r\'eunion dans $\underline{\underline V} (A)$) et l'on note
$\hbox{\LARGE $a$}_X$ la sous-alg\`ebre du produit des
$\hbox{\LARGE $a$}_{X_i}$, $i \in I$, form\'ee des familles
$(\a_i)_{i \in I}$ telles que, si $i \ne j$, les images de $\a_i$
et $\a_j$ dans $\hbox{\LARGE $a$}_{X_{ij}}$ co{\"\i}ncident.

Alors l'objet $X =(\underline{\underline X} , \hbox{\LARGE
$a$}_X)$ sur ${\mathbb F}_1$ (avec les morphismes d'\'evaluation
\break \'evidents) est une vari\'et\'e sur ${\mathbb F}_1$ telle
que
$$
X \build\ot_{{\mathbb F}_1}^{} {\mathbb Z}= V \, .
$$
}

\medskip

\noindent {\bf Preuve.} Les compos\'ees des morphismes d'alg\`ebres
$$
{\mathcal O} (V_{\mathbb C}) \ra \hbox{\LARGE $a$}_{X_i} \ra
\hbox{\LARGE $a$}_{X_{ij}}
$$
et
$$
{\mathcal O} (V_{\mathbb C}) \ra \hbox{\LARGE $a$}_{X_j} \ra
\hbox{\LARGE $a$}_{X_{ij}}
$$
co{\"\i}ncident. Il existe donc une immersion
$$
u : X \ra V
$$
dans ${\mathcal O}$. Par ailleurs, si $W \in {\rm Ob} ({\mathcal
V}_{\mathbb Z})$ et si
$$
\vp : X \ra W
$$
est un morphisme de ${\mathcal O}$, la restriction $\vp_i$ (resp.
$\vp_{ij}$) de $\vp$ \`a $X_i$ (resp. $X_{ij}$) est induite par un
unique morphisme alg\'ebrique $\vp_{i{\mathbb Z}}$ (resp.
$\vp_{ij{\mathbb Z}}$) de $U_i$ (resp. $U_i \cap U_j$) vers $W$.
Comme la restriction de $\vp_{i{\mathbb Z}}$ \`a $X_{ij}$
co{\"\i}ncide avec celle de $\vp$, la restriction de
$\vp_{i{\mathbb Z}}$ \`a $U_i \cap U_j$ est \'egale \`a
$\vp_{ij{\mathbb Z}}$ (unicit\'e). C'est donc aussi la restriction
de $\vp_{j{\mathbb Z}}$ \`a $U_i \cap U_j$. Par cons\'equent, la
collection des morphismes $\vp_{i{\mathbb Z}}$, $i \in I$,
d\'efinit par recollement un morphisme $\vp_{\mathbb Z} \in {\rm
Hom}_{\mathbb Z} (V,W)$ dont la restriction \`a $U_i$ est \'egale
\`a $\vp_{i{\mathbb Z}}$, quel que soit $i \in I$. Par suite, si $A
\in {\rm Ob} ({\mathcal A})$, l'application
$$
\underline{\underline \vp}_{\mathbb Z} \circ \underline{\underline
u} : \underline{\underline X} (A) \ra \underline{\underline W} (A)
$$
co{\"\i}ncide avec $\underline{\underline \vp}_i$ sur le
sous-ensemble $\underline{\underline X}_i (A)$. Comme
$\underline{\underline X} (A)$ est, par d\'efinition, la r\'eunion
des $\underline{\underline X}_i (A)$, $i \in I$, on voit que
$$
\underline{\underline \vp}= \underline{\underline \vp}_{\mathbb Z}
\circ \underline{\underline u} \, .
$$
Par ailleurs, le morphisme d'alg\`ebres compos\'e
$$
\begin{CD}
{\mathcal O} (W_{\mathbb C}) @>{\vp_{\mathbb Z}^*}>>  {\mathcal O}
(V_{\mathbb C})@>{u^*}>> \hbox{\LARGE $a$}_X @>{ \ }>> \hbox{\LARGE
$a$}_{X_i}
\end{CD}
$$
co{\"\i}ncide avec $\vp_i^*$ et $\hbox{\LARGE $a$}_X$ est contenue
dans le produit des alg\`ebres $\hbox{\LARGE $a$}_{X_i}$, donc
$$
\vp^*= u^* \circ \vp_{\mathbb Z}^* \, .
$$
Par cons\'equent $\vp= \vp_{\mathbb Z} \circ u$.

\hfill q.e.d.

\bigskip

\section{Exemples}\label{sec5}

\subsection{ \ } Soient $N \simeq {\mathbb Z}^d$ un ${\mathbb
Z}$-module libre de rang $d \geq 1$, et $\D$ un \'eventail de
$N_{\mathbb R}= N \build\ot_{\mathbb Z}^{} {\mathbb R}$. Notons
${\mathbb P} (\D)$ la vari\'et\'e torique associ\'ee \`a $\D$
\cite{Fu} \cite{Mai}. Supposons que $\D$ est r\'egulier et
que, par cons\'equent, ${\mathbb P} (\D)$ est une vari\'et\'e lisse
sur ${\mathbb Z}$. On se propose de montrer que ${\mathbb P} (\D)$
est l'extension \`a ${\mathbb Z}$ d'une vari\'et\'e $X(\D)$ sur
${\mathbb F}_1$.

Par d\'efinition, $\D$ est une famille finie $\{ \tau \}$ de
c\^ones stricts de $N_{\mathbb R}$. A chaque c\^one $\tau$ est
associ\'ee une vari\'et\'e torique affine
$$
U_{\tau} ={\rm Spec} ({\mathbb Z} \, [{\mathcal S}_{\tau}]) \, ,
$$
o\`u ${\mathcal S}_{\tau}$ est le mono{\"\i}de intersection de $M =
{\rm Hom}_{\mathbb Z} (N , {\mathbb Z})$ avec le dual $\tau^* \sbs
M_{\mathbb R}$ du c\^one $\tau$. La vari\'et\'e ${\mathbb P} (\D)$
est obtenue par recollement des vari\'et\'es $U_{\tau}$ le long de
leurs ouverts $U_{\tau \cap \tau'}$.

Si $m \in {\mathcal S}_{\tau}$ on note
$$
\chi^m : U_{\tau} \ra {\mathbb A}^1
$$
le caract\`ere associ\'e \`a $m$. Si $R \in {\rm Ob} ({\mathcal
R})$, on note $\mu (R)$ l'ensemble (fini) des racines de l'unit\'e
de $R$ et
$$
\underline X_{\tau} (R) \sbs U_{\tau} (R)
$$
l'ensemble (fini) des \'el\'ements $x \in U_{\tau} (R)$ tels que
$$
\chi^m (x) \in \mu (R) \cup \{ 0 \}
$$
quel que soit $m \in {\mathcal S}_{\tau}$.

Par ailleurs, suivant Batyrev et Tschinkel \cite{BT}, on d\'esigne
par $C_{\tau}$ l'ensemble des points $x \in U_{\tau} ({\mathbb C})$
tels que
$$
\vert \chi^m (x) \vert \leq 1
$$
quel que soit $m \in {\mathcal S}_{\tau}$, et
$$
C_{\D} = \, \build\cup_{\tau}^{} C_{\tau}
$$
dans ${\mathbb P} (\D) ({\mathbb C})$. Notons $\hbox{\LARGE
$a$}_{\tau}$ l'alg\`ebre des fonctions complexes continues sur
$C_{\tau}$ dont la restriction \`a l'int\'erieur  de
$C_{\tau}$ est holomorphe. On note aussi $\hbox{\LARGE $a$}_{\D}$
l'alg\`ebre des fonctions complexes continues sur $C_{\D}$ dont la
restriction \`a l'int\'erieur  de
chaque  $C_{\tau}$, $\tau \in \D$, est
holomorphe.

Si $R \in {\rm Ob} ({\mathcal R})$ et si $\s \in R \ra {\mathbb C}$
est un morphisme d'anneaux unitaires, l'application
$$
U_{\tau} (R) \ra U_{\tau} ({\mathbb C})
$$
induite par $\s$ envoie $\underline X_{\tau} (R)$ dans $C_{\tau}$.
On peut donc \'evaluer les fonctions de $\hbox{\LARGE $a$}_{\tau}$
en chaque point $x \in \underline X_{\tau} (R)$. On obtient ainsi
un truc $X_{\tau}= (\underline X_{\tau} , \hbox{\LARGE
$a$}_{\tau})$ sur ${\mathbb F}_1$. Si $A$ est une vari\'et\'e
affine sur ${\mathbb F}_1$ on pose
$$
\underline{\underline{X(\D)}} (A) = \ \build\cup_{\tau}^{} {\rm
Hom}_{\mathbb Z} (A_{\mathbb Z} , U_{\tau})
$$
dans ${\rm Hom}_{\mathbb Z} (A_{\mathbb Z} , {\mathbb P} (\D))$.
Par image inverse et l'\'enonc\'e i) ci-dessous, un point $x \in
\underline{\underline{X(\D)}} (A)$ d\'efinit une \'evaluation $e_x
: \hbox{\LARGE $a$}_{\D} \ra \hbox{\LARGE $a$}_A$.

\medskip

\noindent {\bf Th\'eor\`eme 1.} i) {\it Pour tout c\^one ouvert
$\tau$ de $\D$, le truc $X_{\tau}$ est une vari\'et\'e affine sur
${\mathbb F}_1$ et
$$
U_{\tau} = X_{\tau} \build\ot_{{\mathbb F}_1}^{} {\mathbb Z} \, .
$$
}

\noindent ii) {\it L'objet $X(\D) = (\underline{\underline{X(\D)}}
, \hbox{\LARGE $a$}_{\D})$ sur ${\mathbb F}_1$ est une vari\'et\'e
sur ${\mathbb F}_1$ telle que
$$
X(\D) \build\ot_{{\mathbb F}_1}^{} {\mathbb Z}= {\mathbb P} (\D)
\, .
$$
}

\medskip

\noindent {\bf Preuve.} Montrons d'abord que i) implique ii). Il
r\'esulte de i) que, si $A$ est une vari\'et\'e affine sur
${\mathbb F}_1$ et $\tau \in \D$,
\begin{equation}
\label{eq13}
{\rm Hom}_{\mathcal T} (A , X_{\tau})= {\rm Hom}_{\mathbb Z}
(A_{\mathbb Z} , U_{\tau}) \, .
\end{equation}
Par cons\'equent, la Proposition~5 permet de recoller les
vari\'et\'es $X_{\tau}$ le long des sous-vari\'et\'es $X_{\tau \cap
\tau'}$ pour obtenir une vari\'et\'e sur ${\mathbb F}_1$ dont
l'extension des scalaires \`a ${\mathbb Z}$ est la vari\'et\'e
${\mathbb P} (\D)$. Si $\tau$ et $\tau'$ sont deux c\^ones de $\D$,
les morphismes d'alg\`ebres
$$
\hbox{\LARGE $a$}_{\D} \ra \hbox{\LARGE $a$}_{\tau} \ra
\hbox{\LARGE $a$}_{\tau \cap \tau'}
$$
et
$$
\hbox{\LARGE $a$}_{\D} \ra \hbox{\LARGE $a$}_{\tau'} \ra
\hbox{\LARGE $a$}_{\tau \cap \tau'}
$$
co{\"\i}ncident, donc, d'apr\`es la Proposition~4, et la
d\'efinition de $\underline{\underline{X(\D)}}$, l'objet $X(\D) =
(\underline{\underline{X(\D)}} , \hbox{\LARGE $a$}_{\D})$, muni des
\'evaluations d\'eduites de (\ref{eq13}), est une vari\'et\'e sur
${\mathbb F}_1$ telle que $X(\D) \build\ot_{{\mathbb F}_1}^{}
{\mathbb Z} ={\mathbb P} (\D)$.

Il reste \`a d\'emontrer i). Soit donc $V$ une vari\'et\'e affine
sur ${\mathbb Z}$ et
$$
\vp : X_{\tau} \ra V
$$
un morphisme de ${\mathcal T}$. Supposons d'abord que $V$ soit la
droite affine ${\mathbb A}_{\mathbb Z}^1 = {\rm Spec} ({\mathbb Z}
\, [T])$. On a ${\mathcal O} ({\mathbb A}_{\mathbb C}^1)= {\mathbb
C} \, [T]$ et l'on peut d\'ecrire comme suit la fonction
$$
\a = \vp^* (T) \in \hbox{\LARGE $a$}_{\tau} \, .
$$
Puisque $\D$ est r\'egulier, il existe une base $\{ m_1 , \ldots ,
m_d \}$ de $M$ telle que $\{ m_1 , \ldots ,$ $m_r \}$ soit une
famille g\'en\'eratrice du semi-groupe ${\mathcal S}_{\tau}$, o\`u
$r$ est la dimension de $\tau$. Soit $\tau'$ le c\^one ouvert
${\mathbb R}_+ m_1 + \ldots + {\mathbb R}_+ m_d$. La carte affine
$$
\chi : U_{\tau'} ({\mathbb C}) \ra {\mathbb C}^d
$$
donn\'ee par
$$
\chi (x)= (\chi^{m_1} (x) , \ldots , \chi^{m_d} (x))
$$
identifie $C_{\tau}$ (resp. son int\'erieur)
avec l'ensemble des
points $(x_i)$ de ${\mathbb C}^d$ tels que
$$
\vert x_1 \vert \leq 1 , \ldots , \vert x_{d-r} \vert \leq 1 , \,
\vert x_{d-r+1} \vert= \ldots= \vert x_d \vert =1
$$
(resp.
$$
\vert x_1 \vert < 1 , \ldots , \vert x_{d-r} \vert < 1 , \, \vert
x_{d-r+1} \vert =\ldots= \vert x_d \vert =1) \, .
$$
Autrement dit $C_{\tau}$ (resp. son int\'erieur) est le produit de
$d-r$ disques ferm\'es (resp. ouverts) et de $r$ cercles (cf.
\cite{Mai}, Proposition~3.2.9).

La restriction de $\a (x_1 , \ldots , x_d)$
 au produit des $d$
cercles $\vert x_j \vert =1$, 
$j = 1,\ldots , d$, admet un
d\'eveloppement de Fourier convergent
$$
\a (e^{i \t_1} , \ldots , e^{i \t_d}) = \sum_{I \sbs {\mathbb Z}^d}
a_I \, e^{i I \cdot \t} \, ,
$$
o\`u $a_I \in {\mathbb C}$ et $I \cdot \t \,= \build\sum_{j =
1}^{d} i_{j} \, \t_{j}$ si $I = (i_1 , \ldots , i_d)$. Comme $\a$
est holomorphe dans $C_{\tau}^0$, les coefficients $a_I$ sont nuls
si l'un des indices $i_1 , \ldots , i_{d-r}$ est strictement
n\'egatif.

Par ailleurs, pour tout entier $n \geq 1$, consid\'erons l'anneau
$$
R_{d,n} ={\mathbb Z} \, [T_1 , \ldots , T_d] / (T_j^n=1 , \, 
j = 1 , \ldots , d) \, ,
$$
et le point $x_n \in X_{\tau} (R_{d,n})$ de coordonn\'ees
$$
\chi^{m_j} (x_n)= T_j \in R_{d,n} \, ,
$$
quel que soit $j = 1 , \ldots , d$. Son image
 $\underline{\vp} ( x_n )$
dans ${\mathbb A}_{\mathbb Z}^1 (R_{d,n})= R_{d,n}$ est la classe
d'un polyn\^ome de Laurent
$$
P_n (T_1 , \ldots , T_d)= \sum_{I \in {\mathbb Z}^d} a_I (n)
T_1^{i_1} \ldots T_d^{i_d} \, ,
$$
o\`u les coefficients $a_I (n)$ sont entiers et presque tous nuls.

Tout morphisme $\s : R_{d,n} \ra {\mathbb C}$ est obtenu en
envoyant chaque $T_j$, $j \in J$,
sur une racine $n$-i\`eme $\z_j$ de l'unit\'e. On a
donc, d'apr\`es la condition (\ref{eq5}),
$$
P_n (\z_1 , \ldots , \z_d) = \a (\z_1 , \ldots , \z_d)
$$
d\`es que $\z_1^n =\z_2^n= \ldots =\z_d^n= 1$. Pour tout
multi-indice $I \sbs {\mathbb Z}^d$ le coefficient de Fourier $a_I$
se calcule par la formule
$$
a_I = (2 \pi)^{-d}
\int_{(S^1)^d} \a (e^{i \t_1} , \ldots , e^{i \t_d}) \,
e^{-i I \cdot \t} \, d\t_1 \ldots d \t_d \, .
$$
C'est donc la limite quand $n$ tend vers l'infini de
$$a_I (n) =
n^{-d} \sum_{\z_1^n = \ldots = \z_d^n = 1} P_n (\z_1 , \ldots , \z_n)
\prod_{j =1}^{d} \z_{j}^{-i_{j}}  \, .
$$
Puisque $a_I (n)$ est un entier quel que soit $n \geq 1$, $a_I$
doit \^etre un entier, nul pour presque tout $I$ (puisque la
s\'erie de Fourier d\'efinissant $\a$ converge). Il en r\'esulte
que la fonction $\a =\vp^* (T)$ est un polyn\^ome
$$
P (T_1 , \ldots , T_d) \in {\mathbb Z} \, [T_1 , \ldots , T_{d-r} ,
T_{d-r+1}^{\pm 1} , \ldots , T_d^{\pm 1}] \, ,
$$
c'est-\`a-dire un \'el\'ement de ${\mathbb Z} [{\mathcal
S}_{\tau}]$. Elle d\'efinit donc un morphisme
$$
\vp_{\mathbb Z} : U_{\tau} \ra {\mathbb A}_{\mathbb Z}^1
$$
tel que $\vp^* : {\mathcal O} ({\mathbb A}_{\mathbb C}^1) \ra
\hbox{\LARGE $a$}_{\tau}$ soit le compos\'e de $\vp_{\mathbb Z}^*$
et de la restriction \`a $C_{\tau}$. Si $R \in {\rm Ob} ({\mathcal
R})$ et si $\s : R \ra {\mathbb C}$ est un morphisme d'anneaux
unitaires, il r\'esulte alors de (\ref{eq7}) que le morphisme
compos\'e
$$
\begin{CD}
\underline{X_{\tau}} (R) @>{\underline\vp}>> \underline{\mathbb
A}_{\mathbb Z}^1 (R)= R @>{\s}>> {\mathbb C}
\end{CD}
$$
co{\"\i}ncide avec
$$
\begin{CD}
\underline{X_{\tau}} (R) @>{ \ }>> U_{\tau} (R) @>{\vp_{\mathbb
Z}}>> R @>{\s}>> {\mathbb C} \, .
\end{CD}
$$
Comme le produit des morphismes $\s$ est injectif on voit que $\vp$
est la restriction de $\vp_{\mathbb Z}$ \`a $X_{\tau}$.

Soit maintenant
$$
\vp : X_{\tau} \ra V
$$
un morphisme de ${\mathcal T}$ o\`u $V \sbs {\mathbb A}_{\mathbb
Z}^n$ est une vari\'et\'e affine arbitraire. Soit $W$
la r\'eunion des composantes horizontales de $V$,
c'est \`a dire le spectre de l'anneau des
fonctions de $V$ modulo torsion,
et $W \ra V$ l'inclusion canonique.
Pour tout $R \in {\rm Ob} ({\mathcal
R})$, l'application $W(R) \ra V(R)$ est bijective.
On peut donc supposer que $V=W$, c'est 
\`a dire que la vari\'et\'e
$V$ est plate sur
$\mathbb Z$. Chacune des
coordonn\'ees
$$
\pi_j : V \ra {\mathbb A}_{\mathbb Z}^1 \, ,
$$
$j =1 , \ldots , n$, d\'efinit un morphisme $\pi_j \circ \vp$ de
$X_{\tau}$ vers ${\mathbb A}_{\mathbb Z}^1$ dans ${\mathcal T}$,
dont on vient de voir qu'il est induit par un morphisme $\psi_{j
{\mathbb Z}} \in {\rm Hom}_{\mathbb Z} (U_{\tau} , {\mathbb
A}_{\mathbb Z}^1)$. Notons
$$
\psi_{\mathbb Z} : U_{\tau} \ra {\mathbb A}^n
$$
le produit des morphismes $\psi_{j {\mathbb Z}}$, $j= 1 , \ldots ,
n$. La restriction de $\psi_{\mathbb Z}$ \`a $X_{\tau}$
co{\"\i}ncide avec le morphisme compos\'e
$$
X_{\tau} \build\longra_{}^{\vp} V \longra {\mathbb A}^n \, .
$$
Par ailleurs l'image de $\psi_{\mathbb Z}$ est contenue dans $V$,
car $C_{\tau}$ est Zariski dense dans $U_{\tau} ({\mathbb C})$ et
donc
$$
\psi_{\mathbb Z} (U_{\tau} ({\mathbb C})) \sbs
\overline{\psi_{\tau} (C_{\tau})} \sbs V ({\mathbb C}) \, .
$$
La deuxi\`eme inclusion ci-dessus est due au fait que si $\a \in
{\mathcal O} ({\mathbb A}_{\mathbb C}^d)$, la restriction de $\a
\circ \psi_{\mathbb Z}$ \`a $C_{\tau}$ est \'egale \`a $\vp^* (\a |
V({\mathbb C}))$; elle est donc nulle si $\a$ s'annule sur
$V({\mathbb C})$. Puisque $V$ est plate sur
$\mathbb Z$ on en conclut que $\psi_{\mathbb Z}$ se factorise
par un morphisme alg\'ebrique $\vp_{\mathbb Z} : U_{\tau} \ra V$,
dont la restriction \`a $X_{\tau}$ est \'egale \`a $\vp$. Cela
d\'emontre le th\'eor\`eme.

\subsection{ \ } En prenant pour $\D$ les \'eventails habituels on
obtient les exemples suivants de vari\'et\'es $X(\D)$ sur ${\mathbb
F}_1$.

\subsubsection{ \ } La {\it droite affine} ${\mathbb A}^1$ sur
${\mathbb F}_1$ est la vari\'et\'e affine sur ${\mathbb F}_1$
d\'efinie par
$$
\underline{\mathbb A}^1 (R)= \mu (R) \cup \{ 0 \} \, ,  \
\hbox{si} \ R \in {\rm Ob} ({\mathcal R}) \, ,
$$
(rappelons que $\mu (R)$
d\'esigne l'ensemble des racines de l'unit\'e 
de $R$),
$\hbox{\LARGE $a$}_{{\mathbb A}^1}$ \'etant l'alg\`ebre des
fonctions continues sur le disque unit\'e 
$D = \{ z \in {\mathbb C}
/ \vert z \vert \leq 1 \}$ qui sont holomorphes dans 
l'int\'erieur de $D$, avec
les \'evaluations \'evidentes. On a
$$
{\mathbb A}^1 \build\ot_{{\mathbb F}_1}^{} {\mathbb Z}= {\mathbb
A}_{\mathbb Z}^1 \, .
$$

\subsubsection{ \ } Le {\it groupe multiplicatif} ${\mathbb G}_m$
sur ${\mathbb F}_1$ est la vari\'et\'e affine sur ${\mathbb F}_1$
d\'efinie par
$$
\underline{\mathbb G}_m (R)= \mu (R) \quad \hbox{si} \quad R \in {\rm
Ob} ({\mathcal R}) \, ,
$$
l'alg\`ebre $\hbox{\LARGE $a$}_{{\mathbb G}_m}$ \'etant celle des
fonctions continues sur le cercle unit\'e. Son extension des
scalaires
$$
{\mathbb G}_m \build\ot_{{\mathbb F}_1}^{} {\mathbb Z}=
 {\mathbb
G}_{m,{\mathbb Z}}
$$
est le sch\'ema en groupe multiplicatif sur ${\rm Spec} ({\mathbb
Z})$.

\subsubsection{ \ } On d\'efinit de m\^eme les produits ${\mathbb
A}^a \ts {\mathbb G}_m^b$ sur ${\mathbb F}_1$, $a,b \in {\mathbb
N}$.

\subsubsection{ \ } L'espace projectif ${\mathbb P}^d$ sur
${\mathbb F}_1$ est une vari\'et\'e sur ${\mathbb F}_1$ obtenue par
recollement de $d+1$ espaces affines ${\mathbb A}^d$ sur ${\mathbb
F}_1$. Son extension \`a ${\mathbb Z}$ est l'espace projectif
${\mathbb P}_{\mathbb Z}^d$ de dimension $d$ sur ${\rm Spec}
({\mathbb Z})$. Si $R \in {\rm Ob} ({\mathcal R})$, l'ensemble fini
$\underline{\mathbb P}^d ({\rm Spec} (R))$ est form\'e des points
de ${\mathbb P}_{\mathbb Z}^d (R)$ dont on peut choisir un
syst\`eme de coordonn\'ees homog\`enes dans $(\mu (R) \cup \{ 0
\})^{d+1}$. On a $\hbox{\LARGE $a$}_{{\mathbb P}^d}= {\mathbb C}$ et,
plus g\'en\'eralement, $\hbox{\LARGE $a$}_{\D}= {\mathbb C}$ quand
l'\'eventail $\D$ est complet.

\subsection{ \ } Soit $\Lb$ un ${\mathbb Z}$-module libre de type
fini et $\Vert \cdot \Vert$ une norme hermitienne sur $\Lb
\build\ot_{\mathbb Z}^{} {\mathbb C}$. On pose $\overline{\Lb} =
(\Lb , \Vert \cdot \Vert)$.

Soit $B = \{ x \in \Lb \build\ot_{\mathbb Z}^{} {\mathbb C} / \Vert
x \Vert \leq 1 \}$ la boule unit\'e et $\Phi$ une partie de $(B
\cap \Lb) - \{ 0 \}$ telle que si $v \in (B \cap \Lb) - \{ 0 \}$,
un et un seul des vecteurs $v$ et $-v$ est dans $\Phi$. Si $R \in
{\rm Ob} ({\mathcal R})$ on d\'esigne par $\underline X (R)$
l'ensemble fini des \'el\'ements de $\Lb \build\ot_{\mathbb Z}^{}
R$ qui peuvent s'\'ecrire
\begin{equation}
\label{eq14}
x = \sum_{v \in \Phi} v \ot \z_v \, ,
\end{equation}
o\`u $\z_v \in \mu (R) \cup \{ 0 \}$. L'ensemble fini $\underline X
(R)$ ne d\'epend pas de $\Phi$ et d\'efinit un foncteur
$$
\underline X : {\mathcal R} \ra {\mathcal E}ns \, .
$$

Par ailleurs, notons $\Lb_0 \sbs \Lb$ le r\'eseau engendr\'e par
$\Lb$ et par $C$ le sous-ensemble de $\Lb_0 \build\ot_{\mathbb Z}^{}
{\mathbb C}$ form\'e des vecteurs de norme au plus \'egale \`a $t =
{\rm card} (\Phi)$. Soit $\hbox{\LARGE $a$}$ l'alg\`ebre des
fonctions complexes continues sur $C$ qui sont holomorphes dans
l'int\'erieur de $C$. Si $x \in \underline X (R)$ et si $\s : R \ra
{\mathbb C}$ est un morphisme d'anneaux unitaires, il r\'esulte de
(\ref{eq14}) que
$$
\Vert \s (x) \Vert = \left\Vert \sum_{v \in \Phi} v \ot \s (\z_v)
\right\Vert \leq t \, .
$$
On peut donc \'evaluer les fonctions de $\hbox{\LARGE $a$}$ en $\s
(x)$.

\medskip

\noindent {\bf Proposition 6.} {\it Le couple $X (\overline \Lb) =
(\underline X , \hbox{\LARGE $a$})$, muni des \'evaluations
pr\'ec\'edentes, est une vari\'et\'e affine sur ${\mathbb F}_1$ dont
l'extension des scalaires \`a ${\mathbb Z}$ est le spectre
$X_{\mathbb Z}$ de l'alg\`ebre sym\'etrique du dual de $\Lb_0$.}

\medskip

\noindent {\bf Preuve.} Consid\'erons le morphisme de
groupes ab\'eliens
$$
\pi : {\mathbb Z}^{\Phi} \ra \Lb
$$
tel que
$$
\pi ((n_v))= \sum_{v \in \Phi} n_v \, v \in \Lb \, .
$$
Si ${\mathbb A}^{\Phi}$ est l'espace affine sur ${\mathbb F}_1$ de
rang $t$, $\pi$ induit un morphisme de trucs sur ${\mathbb F}_1$ de
${\mathbb A}^{\Phi}$ vers $X(\overline \Lb)$, qu'on note aussi $\pi$
(on remarquera que si $(z_v) \in {\mathbb C}^{\Phi}$ et $\vert z_v
\vert \leq 1$ quel que soit $v \in \Phi$, le vecteur $\build\sum_{v
\in \Phi}^{} z_v \, v$ est dans $C$).

Si $V \in {\rm Ob} ({\mathcal V}_{\mathbb Z})$ est une vari\'et\'e
affine et si
$$
\vp : X(\overline \Lb) \ra V
$$
est un morphisme de ${\mathcal T}$, le morphisme compos\'e
$$
\vp \circ \pi : {\mathbb A}^{\Phi} \ra V
$$
est, d'apr\`es le Th\'eor\`eme 1, la restriction d'un morphisme
$$
\psi_{\mathbb Z} \in {\rm Hom}_{\mathbb Z} ({\mathbb A}_{\mathbb
Z}^{\Phi} , V) \, .
$$

Le choix d'une section $s$ de la projection ${\mathbb Z}$-lin\'eaire
$\pi$ de ${\mathbb Z}^{\Phi}$ sur $\Lb_0$ induit une section
alg\'ebrique
$$
s_{\mathbb Z} : X_{\mathbb Z} \ra {\mathbb A}_{\mathbb Z}^{\Phi}
$$
de la projection $\pi_{\mathbb Z} \in {\rm Hom}_{\mathbb Z} ({\mathbb
A}_{\mathbb Z}^{\Phi} , X_{\mathbb Z})$. Notons
$$
\vp_{\mathbb Z} : X_{\mathbb Z} \ra V
$$
le morphisme compos\'e $\psi_{\mathbb Z} \circ s_{\mathbb Z}$. Si $\a
\in {\mathcal O} (V_{\mathbb C})$, la restriction de $\vp_{\mathbb
Z}^* (\a)$ \`a $C$ co{\"\i}ncide avec
$$
s^* \circ \pi^* \circ \vp^* (\a) = \vp^* (\a) \, .
$$
Il en r\'esulte, comme dans la fin de la preuve du Th\'eor\`eme 1,
que la restriction de $\vp_{\mathbb Z}$ \`a $X(\overline \Lb)$ est
\'egale \`a $\vp$. \hfill q.e.d.

\bigskip

\noindent {\bf Remarque.} Du point de vue de la th\'eorie d'Arakelov,
le couple $\overline\Lb = (\Lb , \Vert \cdot \Vert)$ est un fibr\'e
vectoriel sur la courbe ${\rm Spec} ({\mathbb Z}) \cup \{ \infty \}$
et les \'el\'ements de $B \cap \Lb$ sont les sections
globales de ce fibr\'e.
C'est un espace vectoriel sur ${\mathbb F}_1$ au sens de Kapranov et
Smirnov (cf. \cite{KS} et 1.3). On peut sans doute voir $X(\overline
\Lb)$ comme la vari\'et\'e affine sur ${\mathbb F}_1$ associ\'e \`a
cet espace vectoriel.

\subsection{ \ } Il faudrait bien s\^ur trouver d'autres exemples que
les pr\'ec\'edents. Par exemple, si $G$ est un sch\'ema en groupes de
Chevalley sur ${\mathbb Z}$, peut-on le d\'efinir sur ${\mathbb
F}_1$~? Et qu'en est-il des vari\'et\'es de drapeaux associ\'ees~?

Par exemple la vari\'et\'e de Grassmann $G(2,4)$ est aussi la conique
$Q$ de ${\mathbb P}_{\mathbb Z}^5$ d'\'equation homog\`ene
\begin{equation}
\label{eq15}
xy - zt + uv = 0 \, .
\end{equation}
Si l'on veut que le nombre des points de $\underline{\underline X}
(R_n)$, $n \geq 1$, v\'erifie le Th\'eor\`eme~2 iii) ci-dessous, on
peut d\'efinir comme suit un objet $X$ sur ${\mathbb F}_1$ contenu
dans celui associ\'e \`a $Q$. Notons $S_{1,{\mathbb Z}}$,
$S_{2,{\mathbb Z}}$, $S_{3,{\mathbb Z}}$ et $S_{4,{\mathbb Z}}$ les
sous-vari\'et\'es localement ferm\'ees de $Q$ suivantes
$$
S_{1,{\mathbb Z}}= Q \cap \{ x \ne 0 \} \simeq {\mathbb A}_{\mathbb
Z}^4
$$
(coordonn\'ees $z/x$, $t/x$, $u/x$, $v/x$),
$$
S_{2,{\mathbb Z}} = Q \cap \{ x= 0 , \, z \ne 0 \} \simeq {\mathbb
A}_{\mathbb Z}^3
$$
(coordonn\'ees $y/z$, $u/z$, $v/z$),
$$
S_{3,{\mathbb Z}} = Q \cap \{ x = z = 0 , \, u \ne 0 \} \simeq
{\mathbb A}_{\mathbb Z}^2
$$
(coordonn\'ees $y/u$, $t/u$), et
$$
S_{4,{\mathbb Z}} = Q \cap \{ x = z = u = 0 \} \simeq 
{\mathbb
P}_{\mathbb Z}^2
$$
(coordonn\'ees homog\`enes $(y,t,v)$).

Notons $S_1 , \, S_2 , \, S_3 , \, S_4$ les vari\'et\'es sur
${\mathbb F}_1$ correspondantes (cf. 5.2), munies des immersions
\'evidentes $S_{\a} \ra Q$ dans ${\mathcal O}$. Si $A \in {\mathcal
A}$ on pose
$$
\underline{\underline X} (A) = \underline{\underline S}_1 (A) \cup
\underline{\underline S}_2 (A) \cup \underline{\underline S}_3 (A)
\cup \underline{\underline S}_4 (A)
$$
dans $\underline{\underline Q} (A)$ et $\hbox{\LARGE $a$}_X =
{\mathbb C}$.

\medskip

\noindent {\bf Question 3.} {\it L'objet 
$X = (\underline{\underline
X} , {\mathbb C})$ sur ${\mathbb F}_1$ est-il une vari\'et\'e sur
${\mathbb F}_1$ telle que $X \build\ot_{{\mathbb F}_1}^{} {\mathbb Z}
= Q$~?}

\section{Fonctions z\^eta}\label{sec6}

\subsection{ \ } Soit $X$ une vari\'et\'e sur ${\mathbb F}_1$. On
cherche \`a lui associer une fonction $\z_X (s)$. Si $n \geq 1$,
rappelons que $R_n = {\mathbb Z} \, [T] / (T^n - 1)$ et que
l'ensemble fini $\underline{\underline X} (R_n)$ est aussi celui des
morphismes de ${\rm Spec} ({\mathbb F}_{1^n})$ vers $X$ (cf.
(\ref{eq8}) si $X$ est affine~; le cas g\'en\'eral se montre de la
m\^eme fa\c con gr\^ace \`a la Proposition~3). Puisque la fonction
z\^eta d'une vari\'et\'e alg\'ebrique sur le corps ${\mathbb F}_q$,
$q > 1$, s'obtient \`a partir du nombre de ses points dans les corps
${\mathbb F}_{q^n}$, $n \geq 1$, il est naturel de d\'efinir $\z_X
(s)$ \`a l'aide des nombres entiers $\# \, \underline{\underline X}
(R_n)$. On fera l'hypoth\`ese simplificatrice suivante~:
\begin{itemize}
\item[(Z)] Il existe un polyn\^ome $N(x) \in {\mathbb Z} \, [x]$ tel
que, quel que soit l'entier $n \geq 1$,
$$
\# \, \underline{\underline X} (R_n) = N (2n+1) \, .
$$
\end{itemize}

En imitant A.~Weil, introduisons alors la s\'erie formelle des
variables $q$ et $T$
$$
Z(q,T) = \exp \left( \sum_{r \geq 1}^{} N(q^r) \, T^r / r \right) \,
.
$$
Pour tout nombre r\'eel $s$ on consid\`ere alors la fonction
$Z(q,q^{-s})$ au voisinage de $q=1$.

\medskip

\noindent {\bf Lemme 1.} i) {\it Pour tout $s \in {\mathbb R}$ la
fonction $Z (q,q^{-s})^{-1}$ a un z\'ero d'ordre $\chi = N(1)$ en
$q=1$. On a de plus
$$
\lim_{q \ra 1} Z (q,q^{-s})^{-1} (q-1)^{-\chi} = \z_X (s) \, ,
$$
 o\`u $\z_X (s)$ est la valeur en $s$ d'un polyn\^ome \`a
coefficients entiers.}

\smallskip

\noindent ii) {\it Si $N(x) \,= \build\sum_{i=0}^{d} a_i \, x^i$ on
a
$$
\z_X (s) = \prod_{i=0}^d (s-i)^{a_i} \, .
$$
}

\noindent {\bf Preuve.} Si $N(x)$ est la somme de deux polyn\^omes
$N'(x)$ et $N''(x)$, les fonctions $Z(q,T)$ et $\z_X (s)$ associ\'ees
\`a $N$ sont les produits de celles associ\'ees \`a $N'$ et $N''$.
Par cons\'equent il suffit de traiter le cas o\`u $N(x) = x^d$. On a
alors
$$
Z(q,T) = \exp \left( \sum_{r \geq1} q^{rd} \, T^r / r \right) = (1 -
q^d \, T)^{-1}
$$
et par cons\'equent
$$
Z(q,q^{-s})^{-1} = 1 - q^{d-s} \, .
$$
Quand $q$ tend vers 1 on trouve
$$
Z(q,q^{-s})^{-1} = (s-d)(q-1) + O(1) \, ,
$$
ce qui d\'emontre le lemme.

\subsection{Th\'eor\`eme 2.} i) {\it Si $\D$ est un \'eventail
r\'egulier, la vari\'et\'e $X(\D)$ sur ${\mathbb F}_1$ qui lui
est associ\'ee (Th\'eor\`eme $1$)
 v\'erifie la condition {\rm (Z)}. Le
polyn\^ome $N(x) \in {\mathbb Z} \, [x]$ tel que
$$
\# \, X(\D) ({\mathbb F}_{1^n}) = N(2n+1) \, , \qquad n \geq 1 \, ,
$$
v\'erifie aussi
$$
\# \, {\mathbb P} (\D) ({\mathbb F}_q) = N(q)
$$
pour tout corps fini ${\mathbb F}_q$ d'ordre $q > 1$. De plus $N(1)$
est la caract\'eristique d'Euler-Poincar\'e de ${\mathbb P} (\D)
({\mathbb C})$.}

\smallskip

\noindent ii) {\it Si $\overline\Lb = (\Lb , \Vert \cdot \Vert)$ est
un r\'eseau hermitien, la vari\'et\'e affine $X(\overline\Lb)$ sur
${\mathbb F}_1$ qui lui est associ\'ee (Proposition $6$) v\'erifie
la condition {\rm (Z)}. Le polyn\^ome $N(x)$ correspondant v\'erifie
$N(1) = 1$ et $N(x) - x^t$ est un polyn\^ome de degr\'e au plus $t-1$
(o\`u $t = {\rm card} (\Phi)$, voir {\rm 5.3}).}

\smallskip

\noindent iii) {\it Soit $Q$ la quadratique d'\'equation (\ref{eq15})
sur ${\mathbb Z}$ et
$$
N(x) = x^4 + x^3 + 2x^2 + x + 1 \, .
$$
On a alors
$$
N(2n+1) = \# \, \underline{\underline X} (R_n)
$$
o\`u $\underline{\underline X}$ est le foncteur d\'efini en {\rm
5.4}, et
$$
N(q) = \# \, Q ({\mathbb F}_q)
$$
pour tout corps fini ${\mathbb F}_q$, $q > 1$. De plus $N(1) = 6$
est la caract\'eristique d'Euler-Poincar\'e de $Q({\mathbb C})$.}

\medskip

\noindent {\bf Preuve.} Pour prouver i) il suffit de traiter le cas
d'une vari\'et\'e torique affine de la forme $U_{\tau}$. En effet, il
en r\'esultera en g\'en\'eral que les nombres $\# \, X(\D) ({\mathbb
F}_{1^n})$ et $\# \, {\mathbb P} (\D) ({\mathbb F}_q)$ sont donn\'ees
par les valeurs (en $2n+1$ et $q$) du m\^eme polyn\^ome, car
$\underline{\underline{X(\D)}} (R_n)$ (resp. ${\mathbb P} (\D)
({\mathbb F}_q)$) est la r\'eunion des ensembles $\underline X_{\tau}
(R_n)$ (resp. $U_{\tau} ({\mathbb F}_q)$), amalgam\'es le long des
sous-ensembles $\underline X_{\tau \cap \tau'} (R_n)$ (resp. $U_{\tau
\cap \tau'} ({\mathbb F}_q)$).

Soit dont $\tau$ un c\^one ouvert r\'egulier de dimension $r$ et $\{
m_1 , \ldots , m_d \}$ une base de $M$ telle que $\{m_1 , \ldots ,
m_r \}$ soit une famille g\'en\'eratrice de ${\mathcal S}_{\tau}$.
Comme dans la preuve du Th\'eor\`eme~1 i), si $\tau' = {\mathbb R}_+
m_1 + \ldots + {\mathbb R}_+ m_d$, la carte
$ U_{\tau'} \ra {\mathbb A}_{\mathbb Z}^d
$
donn\'ee par la famille des $\chi^{m_i}$, $1 \leq i \leq d$,
identifie $\underline X_{\tau} (R)$ \`a l'ensemble fini $(\mu (R)
\cup \{ 0 \})^{d-r} \ts \mu (R)^r$. Cette m\^eme carte identifie
$U_{\tau} ({\mathbb F}_q)$ \`a l'ensemble ${\mathbb F}_q^{d-r} \ts
({\mathbb F}_q^*)^r$. Par cons\'equent le premier \'enonc\'e du
th\'eor\`eme est v\'erifi\'e avec
$$
N(x) = x^{d-r} (x-1)^r \, .
$$
Par ailleurs, $U_{\tau} ({\mathbb C})$ a le type d'homotopie de
$(S^1)^r$. Sa caract\'eristique d'Euler-Poincar\'e est donc nulle si
$r > 0$ est \'egale \`a 1 si $r=0$. On a donc bien $N(1) = \chi
(U_{\tau} ({\mathbb C}))$.

Dans le cas ii) l'ensemble $\underline X (R_n)$ est celui des
\'el\'ements de $\Lb \build\ot_{\mathbb Z}^{} R_n$ de la forme
$$
x = \sum_{v \in \Phi} v \ot \z_v \, ,
$$
o\`u $\z_v \in \mu (R_n) \cup \{ 0 \}$ (voir (\ref{eq14})).
L'ensemble $\mu (R_n)$ est form\'e des polyn\^omes $\pm T^i$,
 $i = 1
, \ldots , n$, et la famille $(T^i)$, $i = 1 , \ldots , n$ est une
base de $R_n$ sur ${\mathbb Z}$, donc tout \'el\'ement de $\Lb
\build\ot_{\mathbb Z}^{} R_n$ s'\'ecrit de fa\c con unique sous la
forme
$$
x = \sum_{i=1}^n x_i \ot T^i
$$
avec $x_i \in \Lb$. Ceci conduit \`a d\'ecrire $\underline X (R_n)$
de la fa\c con suivante. Si $\Phi' \sbs \Phi$ est une partie de
$\Phi$, appelons $V(\Phi')$ l'ensemble des vecteurs non nuls de $\Lb$
de la forme
$$
\sum_{v \in \Phi'} \ve_v \, v \, , \quad \hbox{avec} \ \ve_v \pm 1
\, .
$$
Si $k \geq 1$ est un entier on note $V(k)$ l'ensemble des collections 
$\{v_1,\ldots,v_k\}$ de $k$ vecteurs de $\Lb$ telles qu'il existe une 
famille finie $(\Phi_1,\ldots,\Phi_k)$ de sous-ensembles disjoints
de $\Phi$ tels que $v_j\in V(\Phi_j), 1 \leq j \leq k$. 
On note $X(k,n)$ l'ensemble des
 \'el\'ements de $\underline X (R_n)$ de la forme
$$
x = \sum_{j=1}^k v_j \ot T^{n_j} \, ,
$$
o\`u $\{v_1,\cdots,v_k\} \in V (k)$ et
 $1\leq n_1 < n_2 < \ldots < n_k \leq n$.
L'ensemble $\underline X (R_n)$
est la r\'eunion disjointe
de $\{ 0 \}$ et des sous-ensembles $X(k,n)$. Par
ailleurs,
$$
\# \, X(k , n) = (\# \, V (k)) n(n-1) \ldots (n-k+1) \, .
$$
Si $v_j \in \Phi_j$ on a $-v_j \in \Phi_j$,
donc le cardinal de chaque ensemble $V(k)$ est divisible
par $2^k$, et $\# \, X
(k, n)$ est un polyn\^o\-me entier de la variable $2n$, nul \`a
l'origine et de degr\'e $k$. 
Si $k=t = {\rm card} (\Phi)$, chaque
ensemble $\Phi_j$ comporte un seul vecteur et son oppos\'e, donc $\#
\, X(t, n)$ est la somme de $(2n)^t$ et d'un polyn\^ome de degr\'e
au plus $t-1$ en la variable $2n$. Il en r\'esulte que $\# \,
\underline X (R_n) = N(2n+1)$ o\`u $N(x) \in {\mathbb Z} \, [x]$
v\'erifie les conditions de l'\'enonc\'e (on notera que $N(1) = 1$
est encore la caract\'eristique d'Euler-Poincar\'e de $X_{\mathbb Z}
({\mathbb C})$).

Pour montrer iii) il suffit de noter que $\underline{\underline X}
(R_n)$ (resp. $Q({\mathbb F}_q)$) est la r\'eunion disjointe des
ensembles $\underline{\underline S}_{\a} (R_n)$ (resp. $S_{\a ,
{\mathbb Z}} ({\mathbb F}_q)$), $\a = 1,2,3,4$, et d'appliquer le
Th\'eor\`eme~2 i).

\subsection{Exemples.}

\subsubsection{ \ } Le point ${\rm Spec} ({\mathbb F}_1)$ a pour
fonction z\^eta
$$
\z_{{\rm Spec} ({\mathbb F}_1)} (s) = s \, .
$$

\subsubsection{ \ } La droite affine ${\mathbb A}^1$ est telle que
$N(q) = q$ et
$$
\z_{{\mathbb A}^1} (s) = s-1 \, .
$$

\subsubsection{ \ } Le groupe multiplicatif ${\mathbb G}_m$ v\'erifie
$N(q) = q-1$ et
$$
\z_{{\mathbb G}_m} (s) = (s-1) / s \, .
$$

\subsubsection{ \ } Si $\D$ et $\D'$ sont deux \'eventails
on a
$$
\z_{X(\D \ts \D')} (s) = \z_{X(\D)} (s) \ts \z_{X(\D')} (s) \, .
$$

\subsubsection{ \ } L'espace projectif ${\mathbb P}^d$ v\'erifie
$$
N(q) = [d]
$$
(voir 1.1) et donc
$$
\z_{{\mathbb P}^d} (s)= s(s-1) \ldots (s-d) \, .
$$
Cette formule et les pr\'ec\'edentes sont celles pr\'evues par Manin
\cite{Ma} (voir (\ref{eq3})).

\subsubsection{ \ } Soit $\Lb = {\mathbb Z}$ et
 $\lb = \Vert 1 \Vert >
0$ la norme des g\'en\'erateurs de $\Lb$. L'entier $t \geq 0$ est la
partie enti\`ere de $\lb$ et la fonction z\^eta de $X(\overline\Lb)$
ne d\'epend que de $t$, i.e. $\z_{X(\overline\Lb)} (s)= \z_t (s)$,
avec $\z_0 (s)= s$, $\z_1 (s) = s-1$, $\z_2 (s) = s(s-2)/(s-1)$,...
 Je ne
connais pas de formule g\'en\'erale pour $\z_t (s)$.

\subsection{Remarques.}

\subsubsection{ \ } La d\'efinition du Lemme 1 i) a peut-\^etre un sens
dans des situations o\`u (Z) n'est pas v\'erifi\'ee. On peut
aussi se demander si une vari\'et\'e $X$ lisse sur ${\mathbb F}_1$
v\'erifie la condition (Z), si $N(q)$ est le cardinal de $(X
\build\ot_{{\mathbb F}_1}^{} {\mathbb Z}) ({\mathbb F}_q)$, $q > 1$,
et si $N(1)$ est la caract\'eristique d'Euler-Poincar\'e de $(X
\build\ot_{{\mathbb F}_1}^{} {\mathbb Z}) ({\mathbb C})$. En termes de
motifs~:

\medskip

\noindent {\bf Question 4.} {\it Les motifs des vari\'et\'es lisses
sur ${\mathbb F}_1$ sont-ils des motifs de Tate mixtes~?}

\subsubsection{ \ } Quand $X = {\rm Spec} (R)$, $R \in {\rm Ob}
({\mathcal R})$, le cardinal de $\underline X (R_n)$ n'a pas un
comportement simple et je ne sais pas d\'efinir $\z_X (s)$. Mais on
remarquera que $X$ n'est jamais lisse sur ${\mathbb F}_1$ (sauf quand
$R = {\mathbb Z}$).

\subsubsection{ \ } Soit $G$ un graphe fini. M.~Kontsevich associe \`a
$G$ des vari\'et\'es $Y_G$ et $X_G$ sur ${\mathbb Z}$, qui sont
``souvent'' des vari\'et\'es polynomialement d\'enombrables \cite{BB},
c'est-\`a-dire telles que leur nombre de points dans ${\mathbb F}_q$,
$q > 1$, est la valeur en $x=q$ d'un polyn\^ome $N(x)$
\`a coefficients entiers. Kontsevich a
sugg\'er\'e que ces vari\'et\'es sont d\'efinies sur ${\mathbb F}_1$.
Lorsquelles sont polynomialement d\'enombrables peut-\^etre
v\'erifient-elles la condition (Z) avec le m\^eme polyn\^ome $N(x)$.
De m\^eme on peut se demander si un matro{\"\i}de d\'efinit
 une vari\'et\'e
sur ${\mathbb F}_1$.

\section{Sur l'image du $J$-homomorphisme}\label{sec7}

\subsection{ \ } Comme l'a not\'e Manin \cite{Ma}, la formule
(\ref{eq1}) sugg\`ere que les groupes de $K$-th\'eorie alg\'ebrique de
${\mathbb F}_1$ sont les groupes d'homotopie stable des sph\`eres. En
effet, rappelons que, d'apr\`es D.~Quillen, la $K$-th\'eorie d'un
anneau $A$ est
\begin{equation}
\label{eq16}
K_m (A) = \pi_m \, {\rm BGL} (A)^+ \, , \quad \hbox{si} \ m \geq 1 \, ,
\end{equation}
o\`u le CW-complexe ${\rm BGL} (A)^+$ est obtenu en adjoignant au
classifiant du groupe lin\'eaire infini ${\rm GL} (A)$ des cellules de
dimensions 2 et 3 (voir \cite{Lo}). Si $A$ est un corps
et $m > 1$, la formule
(\ref{eq16}) vaut aussi en rempla\c cant le groupe lin\'eaire par le
groupe sp\'ecial lin\'eaire.

Par ailleurs, un th\'eor\`eme de Barratt, Priddy et Quillen \cite{P}
affirme que si $\Si_{\ify} = \build\cup_{N \geq 1}^{} \Si_N$ est le
groupe sym\'etrique infini, on a
\begin{equation}
\label{eq17}
\pi_m \, B \, \Si_{\ify}^+ = \pi_m^s \, , \quad m \geq 1 \, ,
\end{equation}
o\`u $\pi_m^s$ est le $m$-i\`eme groupe d'homotopie stable des
sph\`eres. Il est donc logique d'\'ecrire
\begin{equation}
\label{eq18}
K_m ({\mathbb F}_1) = \pi_m^s \, , \quad m \geq 1 \, .
\end{equation}

\medskip

{\bf Variante.} Le groupe des unit\'es de ${\mathbb F}_1$ est
${\mathbb F}^*_1 = \mu ({\mathbb Z}) = \{ \pm 1 \}$. On peut donc
envisager que ${\rm GL_N} ({\mathbb F}_1)$ est le groupe des matrices
monomiales $N \ts N$ dont les entr\'ees sont $\pm 1$, c'est-\`a-dire
le produit en couronnes $\Si_N \int ({\mathbb Z} / 2)^N$. Cela
conduirait \cite{K}
 \`a remplacer (\ref{eq18}) par la formule
$$
K_m ({\mathbb F}_1) = \pi_m^s (B ({\mathbb Z} / 2)) \, ,
$$
qui n'en diff\`ere que par un groupe fini de 2-torsion.

\subsection{ \ } L'inclusion standard de $\Si_N$ dans ${\rm GL}_N
({\mathbb Z})$ induit un morphisme
$$
\a_m : \pi_m^s \ra K_m ({\mathbb Z}) \, , \quad m \geq 1 \, ,
$$
que l'on peut voir comme celui induit en $K$-th\'eorie alg\'ebrique
par le morphisme d'anneaux ${\mathbb F}_1 \ra {\mathbb Z}$. Ce
morphisme $\a_m$ est bien compris.

Adams a introduit un morphisme
$$
J : \pi_m \, O \ra \pi_m^s \, ,
$$
dont il montre que l'image est un groupe cyclique. Le groupe ${\rm
Im} \,( J )$ est d'ordre 2 si $m$ est congru \`a 0 ou 1 modulo 8. Il
est cyclique d'ordre le d\'enominateur $w_i$ de $b_i / 2i$ si $m =
2i-1$ avec $i$ pair, o\`u $b_i$ le $i$-\`eme nombre de Bernoulli. Et
il est nul sinon.

Quillen a montr\'e que $\a_m$ est injectif sur l'image de $J$ \cite{Q}
et S.A. Mitchell \cite{Mi}
a montr\'e que l'image de $\a_m$ co{\"\i}ncide avec
celle de $\a_m \circ J$.

\subsection{ \ } Par ailleurs, les th\'eorie de M.~L\'evine, M.
Hanamura et V. V\oe vodsky permettent de voir les groupes de
$K$-th\'eorie alg\'ebrique comme des groupes d'extensions dans une
cat\'egorie (d\'eriv\'ee) de motifs mixtes. Si par exemple
 $m=2i-1 >
1$, on a
\begin{equation}
\label{eq19}
K_{2i-1} ({\mathbb Z}) = K_{2i-1} ({\mathbb Q}) = {\rm
Ext}_{{\rm DM}({\mathbb Q})} ({\mathbb Z}(i) , {\mathbb Z}) \, .
\end{equation}
On peut s'attendre \`a une formule du m\^eme type pour le corps \`a un
\'el\'ement~:
$$
K_{2i-1} ({\mathbb F}_1)= {\rm Ext}_{{\rm DM} ({\mathbb F}_1)}
({\mathbb Z} (i), {\mathbb Z}) \, .
$$
Les r\'esultats de 7.2 signifieraient alors que la classe d'une
extension de motifs de Tate sur ${\mathbb Q}$ qui provient, par
extension des scalaires de ${\mathbb F}_1$ \`a ${\mathbb Z}$ (puis
${\mathbb Q}$), d'une extension de motifs de Tate sur ${\mathbb F}_1$
est de torsion, et annul\'ee par $w_i$ s'il s'agit d'une extension de
${\mathbb Z}$ par ${\mathbb Z} (i)$ (ou de ${\mathbb Z} (j)$ par
${\mathbb Z} (i+j)$, $j \in {\mathbb Z}$).

\subsection{ \ } Un r\'esultat de B.~Totaro \cite{To} est coh\'erent
avec les r\'eflexions pr\'ec\'edentes. Si ${\mathbb P} (\D)$ est la
vari\'et\'e torique associ\'ee \`a un \'eventail $\D$ de dimension
$d$, le Th\'eor\`eme 5 de \cite{To} montre que la filtration des poids
de la cohomologie singuli\`ere \`a coefficients rationnels de
${\mathbb P} (\D) ({\mathbb C})$ est canoniquement scind\'ee. La
preuve consiste \`a consid\'erer la filtration canonique
$$
{\mathbb P} (\D)= Y_0 \supset Y_1 \supset \ldots \supset Y_d \supset
\emptyset \, ,
$$
o\`u $Y_k$ est la r\'eunion des orbites toriques de dimension au plus
$d-k$. Chaque strate $Y_k - Y_{k-1}$ est la r\'eunion disjointe de
produits du groupe multiplicatif par la droite affine, et la
multiplication par $n$ dans $N \simeq {\mathbb Z}^d$ (cf.~5.1) induit le
produit par $n^j$ sur la partie de poids $j$ de la cohomologie de ces
strates.

Autrement dit, la filtration de ${\mathbb P} (\D)$ par $Y_k$ fournit
des extensions de motifs de Tate, et donc des classes dans les groupes
${\rm Ext}_{{\rm DM}({\mathbb Q})} ({\mathbb Z} (i+j) , {\mathbb Z}
(j))$, qui sont annul\'ees par le p.g.c.d. des entiers $n^{i+j} - n^j$,
$n > 1$. Si $j$ est assez grand, on sait bien que ce p.g.c.d. n'est
autre que $w_i$. Or le Th\'eor\`eme~1 indique que ces extensions sont
sans doute dans l'image du morphisme
$$
{\rm Ext}_{{\rm DM} ({\mathbb F}_1)} ({\mathbb Z} (i+j) , {\mathbb Z}
(j)) \ra {\rm Ext}_{{\rm DM} ({\mathbb Q})} ({\mathbb Z} (i+j) ,
{\mathbb Z} (j)) \, ,
$$
ce qui va dans le sens de la discussion de 7.3.

On peut aussi se demander si le r\'esultat de Totaro est optimal et si
le groupe ${\rm Im} (\a_{2i-1})$, $i \geq 1$, est engendr\'e par des
extensions de motifs sur ${\mathbb Q}$ provenant des vari\'et\'es
toriques par le proc\'ed\'e pr\'ec\'edent. Ceci conduit au probl\`eme
suivant, que l'on pourrait aborder en \'etudiant la structure de Hodge
mixte \`a coefficients entiers de ${\mathbb P} (\D)$ et de sa
filtration canonique par les $Y_k$~:

\medskip

\noindent {\bf Question 5.} {\it Pour tout entier $i \geq 1$,
existe-t-il $j \geq 0$, une vari\'et\'e torique ${\mathbb P} (\D)$ et
une extension de ${\mathbb Z} (j)$ par ${\mathbb Z} (i+j)$, d\'eduite
de la filtration canonique de ${\mathbb P} (\D)$, dont la classe dans
${\rm Ext}_{{\rm DM} ({\mathbb Q})} ({\mathbb Z} (j) , {\mathbb Z}
(i+j))$ soit d'ordre $w_i$~?}


\newpage

\begin{thebibliography}{99}
\bibitem{A} {\sc E.Arbarello, C. De Concini et V.C.  Kac},
The infinite wedge representation and the reciprocity law for algebraic curves,
Theta functions, Proc. 35th Summer Res. Inst. Bowdoin Coll., Brunswick/ME 1987, 
{\it Proc. Symp. Pure Math.} {\bf 49}, Pt. 1,  (1989) 171-190.
\bibitem{BT} {\sc V. Batyrev et Yu. Tschinkel}, Height zeta functions
of toric varieties, {\it J. Math. Sci.}, New York {\bf 82}, No.1,
(1996) 3220-3239.
\bibitem{BB} {\sc P.Belkale et P.Brosnan}, Matroids, motives and a conjecture
of Kontsevich, {\it Duke Math. J.}, {\bf 116}, No.1,
(2003) 147-188.
\bibitem{B1} {\sc M.~Brou\'e, G.~Malle et J.~Michel}, Th\'eor\`emes de
Sylow g\'en\'eriques pour les groupes r\'eductifs sur les corps finis,
{\it Math. Ann.} {\bf 292}, No.2, (1992) 241-262.
\bibitem{B2} {\sc M.~Brou\'e, G.~Malle et J.~Michel},
Repr\'esentations unipotentes g\'en\'eriques et blocs des groupes
r\'eductifs finis avec un appendice de George Lusztig, {\it
Ast\'erisque} {\bf 212}; Montrouge: {\it Soci\'et\'e
Math\'ematique de France}
203 p. (1993).
\bibitem{B3} {\sc M.~Brou\'e, G.~Malle et J.~Michel}, Towards Spetses
I, {\it Transform. Groups} {\bf 4}, No.2-3, (1999) 157-218.
\bibitem{C} {\sc A. Connes}, Sym\'etries Galoisiennes et
Renormalisation, Pr\'e-publication IHES/M/02/79 (2002).
\bibitem{DG} {\sc M. Demazure et P. Gabriel}, Introduction to
algebraic geometry and algebraic groups, {\it North-Holland
Mathematics Studies} {\bf 39} Amsterdam, New York, Oxford:
North-Holland Publishing Company. XIV, 358 p. (1980).
\bibitem{Fu} {\sc W. Fulton}, {\it Introduction to toric
varieties}, Princeton University Press (1993).
\bibitem{K} {\sc M. Kapranov}, Some conjectures on the 
absolute direct image, lettre, 26/05/1995.
\bibitem{KS} {\sc M. Kapranov et A. Smirnov},
 Cohomology determinants
and reciprocity laws: number field case, pr\'epublication.
\bibitem{Lo} {\sc J.-L. Loday}, K-th\'eorie alg\'ebrique et
repr\'esentations de groupes, {\it Ann. Sci. \'Ec. Norm. Sup\'er.},
IV. S\'er. 9 (1976) 309-377.
\bibitem{Mai} {\sc V. Maillot}, G\'eom\'etrie d'Arakelov des
vari\'et\'es toriques et fibr\'es en droites int\'egrables, {\it
M\'em. Soc. Math. Fr.} Nouv. S\'er. 80, 129 p. (2000).
\bibitem{Ma} {\sc Y. Manin}, Lectures on zeta functions and motives
(according to Deninger and Kuro\-kawa), Columbia University Number
Theory Seminar (New York, 1992); {\it Ast\'erisque} {\bf 228}
 (1995), 4,
121--163.
\bibitem{Mi} {\sc S.A. Mitchell}, The Morava K-theory of algebraic
K-theory spectra,{\it J. K-Theory} {\bf 3}, No.6, (1990) 607-626.
\bibitem{P} {\sc S.B. Priddy}, On $\Omega^\infty S^\infty$ and the
infinite symmetric group, {\it Algebraic Topology, Proc. Sympos. Pure
Math.} {\bf 22}, (1971) 217-220.
\bibitem{Q} {\sc D. Quillen}, Letter from Quillen to Milnor on Im
$(\pi_i0\rightarrow\pi^s_i\rightarrow K_i{\mathbb Z})$, Algebraic
K-Theory, Proc. Conf. Evanston 1976, {\it Lect. Notes Math.} {\bf
551}, (1976) 182-188.
\bibitem{R} {\sc W. Rudin}, Real and complex analysis, 3rd ed.
New York, NY McGraw-Hill, xiv, 416 p.  (1987).
\bibitem{S} {\sc A. Smirnov}, Hurwitz inequalities for number
fields, {\it St. Petersbg. Math. J.} {\bf 4}, No.2, (1993) 357-375;
translation from {\it Algebra Anal.} {\bf 4}, No.2, (1992) 186-209.
\bibitem{So} {\sc C. Soul\'e}, On the field with one element
(expos\'e \`a l'Arbeitstagung, Bonn, June 1999), Preprint IHES/M/99/55
(1999).
\bibitem{St} {\sc R. Steinberg}, A geometric approach to
 the representations of the full linear group over a Galois field,
 {\it Trans. Am. Math. Soc.} {\bf 71}, 274-282 ,(1951)
\bibitem{T} {\sc J. Tits}, Sur les analogues alg\'ebriques des groupes
semi-simples complexes. Colloque d'alg\`ebre sup\'erieure, Bruxelles,
1956, pp. 261-289. Centre Belge de Recherches Math\'ematiques,
Gauthier-Villars (1957).
\bibitem{To} {\sc B. Totaro}, Chow groups, Chow cohomology, and linear
varieties, preprint (1998).
\bibitem{W} {\sc A. Weil}, Sur l'analogie entre les corps de nombres
alg\'ebriques et les corps de fonctions  alg\'ebriques, [1939a] in
{\it Oeuvres Scient.\/} {\bf I}, 236-240, (1980), Springer-Verlag.
\end{thebibliography}
\end{document}
--MIMEStream=_0+71920_64517087713423_58240431625--